\documentclass[12pt,reqno]{amsart}
\usepackage{amsmath}
\usepackage{amssymb}
\usepackage{amsthm}

\newtheorem{theorem}{Theorem}[section]
\newtheorem{prop}[theorem]{Proposition}
\newtheorem{lem}[theorem]{Lemma}
\newtheorem*{cor}{Corollary}
\theoremstyle{definition}
\newtheorem{defn}[theorem]{Definition}
\newtheorem{ex}[theorem]{Example}
\theoremstyle{remark}
\newtheorem*{rem}{Remark}
\numberwithin{equation}{section}

\newcommand{\Gg}{\mathfrak{g}}    
\newcommand{\Gh}{\mathfrak{h}}
\newcommand{\Gq}{\mathfrak{q}}
\newcommand{\Gm}{\mathfrak{m}}

\newcommand{\Gz}{\mathfrak{z}}

\begin{document}

\title[Non-compact quantum groups]
{Non-compact quantum groups arising from Heisenberg type Lie bialgebras}

\author{Byung--Jay Kahng}
\date{}
\address{Department of Mathematics\\ University of California at Davis\\
Davis, CA 95616}
\email{bkahng@math.ucdavis.edu}

\begin{abstract}
The dual Lie bialgebra of a certain {\em quasitriangular\/} Lie bialgebra
structure on the Heisenberg Lie algebra determines a (non-compact)
Poisson--Lie group $G$.  The compatible Poisson bracket on $G$ is
non-linear, but it can still be realized as a ``cocycle perturbation''
of the linear Poisson bracket.  We construct a certain twisted group 
$C^*$--algebra $A$, which is shown to be a strict deformation quantization
of $G$.  Motivated by the data at the Poisson (classical) level, we then
construct on $A$ its {\em locally compact quantum group\/} structures:
comultiplication, counit, antipode and Haar weight, as well as its
associated multiplicative unitary operator.  We also find a quasitriangular
``quantum universal $R$--matrix'' type operator for $A$, which agrees well
with the quasitriangularity at the Lie bialgebra level.
\end{abstract}
\maketitle

{\sc Introduction.} 
 So far, usual method of constructing quantum groups has been the method of
``generators and relations'', in which one tries to {\em deform\/} the relations
between the generators (i.\,e. ``coordinate functions'') of the commutative
algebra of functions on a Lie group.  But if we wish to study locally compact
($C^*$--algebraic) quantum groups, this provides a serious obstacle:  In
non-compact situations, the generators tend to be unbounded, which makes it
difficult to treat them in the $C^*$--algebra framework (For example see
\cite{Wr4}, where Woronowicz introduces the highly technical theory of
``unbounded operators affiliated with $C^*$--algebras'' in his construction
of the quantum $E(2)$ group.).  In addition, the method of generators and
relations is at best an indirect method, in the sense that the deformation of
the pointwise product on the function algebra is not explicitly obtained.

 Because of this, constructing new (especially, non-compact) quantum groups
has been rather difficult.  Among the specific examples of non-compact quantum
groups which have been constructed and studied are: \cite{PW}, \cite{Wr4,Wr6},
\cite{Bj}, \cite{VD}, \cite{SZ}, \cite{Rf5,Rf7}, \cite{Zk2}, \cite{Ld}.

 Recently in \cite{BJKp1}, we defined certain (in general non-linear) Poisson
brackets on dual vector spaces of Lie algebras, denoted by $\{\ ,\ \}_{\omega}$,
which are ``cocycle perturbations'' of the linear Poisson brackets.  We then
showed that deformation quantizations of these Poisson brackets (of which
the linear Poisson brackets are special cases) are provided by twisted
group $C^*$--algebras.  This construction is relatively general (at least
for those Poisson brackets of the aforementioned type), and there are some
indications \cite{BJKp1} that further generalization could be possible.  In
addition, it is a direct approach, where we deform the pointwise product
directly at the function algebra level.

 We wish to use this method to construct some $C^*$--algebraic quantum groups.
But to construct a quantum group from a twisted group $C^*$--algebra, it should
be given a compatible comultiplication and other quantum group structures.  If
we are to reasonably expect a twisted group $C^*$--algebra  (regarded here as
a deformation quantization of our Poisson bracket $\{\ ,\ \}_{\omega}$) to be
also equipped with a  compatible comultiplication, we need to require that
$\{\ ,\ \}_{\omega}$ determines a Poisson--Lie group. 

 Since a typical Poisson bracket we consider is defined on the dual space of
a Lie algebra, this means that it is reasonable to impose a condition that
the dual vector space is itself a Lie group such that it forms, together with
the given Poisson bracket, a Poisson--Lie group.  This suggests us to consider
the following.

 Suppose $H$ is a Poisson--Lie group.  Then its Lie algebra $\Gh$ is a Lie
bialgebra such that its dual vector space $\Gg=\Gh^*$ is also a Lie bialgebra.
The Lie group $G$ of $\Gg$ is the {\em dual Poisson--Lie group\/} of $H$ (See
\cite{LW}, \cite{CP}, or Appendix of \cite{BJK} for discussion on Poisson--Lie
groups.).  In other words, at the level of Poisson--Lie groups, the notion of
a Poisson bracket defined on the ``dual'' of a Lie group naturally exists.
Moreover, if the dual Poisson--Lie group $G$ is exponential solvable (so $G$
is diffeomorphic to $\Gg$ via the exponential map), then we may transfer via
the exponential map the compatible Poisson bracket on $G$ to a Poisson bracket
on $\Gg$.  To apply the result of \cite{BJKp1}, let us assume that the resulting
Poisson bracket on $\Gg=\Gh^*$ be of our type discussed above.  

 Then by the main theorem (Theorem 3.4) in \cite{BJKp1}, a deformation quantization
of $\Gg$ is given in terms of the twisted group $C^*$--algebra of $H$.  Since
$\Gg\cong G$, this can also be regarded as a deformation quantization of the
Poisson--Lie group $G$.  In particular, if $G$ is globally linearizable (i.\,e.
the compatible Poisson bracket on $G$ is Poisson isomorphic to the linear Poisson
bracket on $\Gg$), its deformation quantization is given by the ordinary group
$C^*$--algebra $C^*(H)$.  

 This set-up does not automatically provide a compatible comultiplication on the
twisted group $C^*$--algebra.  But we can usually collect enough data at the
Poisson--Lie group level so that the candidates for comultiplication and other
quantum group structures could be obtained.  We then have to provide a rigorous
analytic proof for our choice of comultiplication, which is not necessarily simple.
It often helps to find some useful tools like multiplicative unitary operators
(in the sense of Baaj and Skandalis \cite{BS}).  

 Many of the earlier known examples of non-compact quantum groups, including
the ones in \cite{VD}, \cite{Rf5}, \cite{SZ}, are deformations of some ``globally
linearizable'' Poisson--Lie groups.  So these quantum groups essentially look
like ordinary group ($C^*$--)algebras (See also \cite[\S7]{Zk2}, \cite{Ld},
\cite{EV}.).  Whereas our method allows us to deform a more general type of
Poisson--Lie groups whose compatible Poisson brackets are in general non-linear.
In fact, these early examples are special cases of our construction.

 In this paper, we will follow the method outlined above to construct some
specific examples of quantum groups.  We will begin the first section with
the study of the $2n+1$ dimensional Heisenberg Lie group $H$, equipped with
a certain (linear) Poisson--Lie group structure (By \cite{SZ}, it is actually
known that all possible Poisson brackets on $H$ are linear).  In particular,
we will consider the one obtained from a certain ``(quasitriangular) classical
$r$--matrix''.

 Then we consider the dual Poisson--Lie group $G$ of $H$.  The dual Poisson
bracket is in general not linear.  But in our case, we show that it is of
the ``cocycle perturbation'' type mentioned earlier.  So following the method
of \cite{BJKp1}, we construct (in section 2) a $C^*$--algebra which is a
deformation quantization of this dual Poisson bracket.

 On this $C^*$--algebra, we construct its quantum group structures, including
comultiplication and the associated multiplicative unitary operator (in section 3),
counit and antipode (in section 4), and Haar weight (in section 5).  In the last
section, we find a {\em quasitriangular\/} ``quantum universal $R$--matrix''
type operator for our $C^*$--algebra, and relate it with the classical
$r$--matrix we started with.

 We discuss the representation theory of our quantum groups in a separate paper
\cite{BJKpz}.  The quantum $R$--matrix operator plays an important role here.
Discussion of more general quantum groups which can be constructed using similar
techniques are also postponed to a future occasion.  For instance, we could
consider a more general two-step nilpotent Lie group whose center has dimension
higher than one, and try to deform its dual Poisson--Lie group equipped with its
non-linear Poisson bracket.  See \cite{BJK} for a discussion.  

 Most of the material in this article formed an integral part of the author's
Ph.D. thesis at U.C. Berkeley \cite{BJK}.  Several modifications and some new
additions were made, including the discussion on the $R$--matrix.  I would 
like to use this opportunity to express my deepest gratitude to my advisor,
Professor Marc Rieffel.  Without his constant encouragement and show of interest,
this work would not have been made possible.  I also thank the referee, who gave
me some corrections and many other helpful suggestions.

\section{The Lie bialgebras, Poisson--Lie groups}

 The notion of Poisson--Lie groups is more or less equivalent to the notion
of Lie bialgebras (\cite{Dr}, \cite{LW}), and these are the objects to be
quantized to produce quantum groups.  In this section, we will study these
``classical'' objects, to find enough data we can use to construct our
specific quantum groups.   The Lie bialgebras we will exclusively study are
either nilpotent or exponential solvable ones, so that we are able to treat
their deformation quantizations in the $C^*$--algebra framework (See \cite
{BJKp1}, \cite{Rf3}.).

\begin{defn} Let $\Gh$ be the $2n+1$ dimensional (real) Lie algebra
generated by $\mathbf{x_i},\mathbf{y_i} (i=1,\dots,n),\mathbf{z}$, 
with the following relations:
$$
[\mathbf{x_i},\mathbf{y_j}]=\delta_{ij}\mathbf{z}, 
\quad [\mathbf{z},\mathbf{x_i}]=[\mathbf{z},\mathbf{y_i}]=0.
$$
This is actually the well-known {\em Heisenberg Lie algebra\/}.  Let us
also consider the extended Heisenberg Lie algebra $\tilde{\Gh}$, generated
by $\mathbf{x_i},\mathbf{y_i} (i=1,\dots,n),\mathbf{z},\mathbf{d}$, with
the relations:
$$
[\mathbf{x_i},\mathbf{y_j}]=\delta_{ij}\mathbf{z},\quad 
[\mathbf{d},\mathbf{x_i}]=\mathbf{x_i},\quad [\mathbf{d},\mathbf{y_i}]
=-\mathbf{y_i},\quad [\mathbf{z},\mathbf{x_i}]=[\mathbf{z},\mathbf{y_i}]
=[\mathbf{z},\mathbf{d}]=0.
$$
\end{defn}

\begin{defn} The (connected and simply connected) Lie group corresponding
to $\Gh$ is the {\em Heisenberg Lie group\/}, denoted by $H$.  The space for
this Lie group is isomorphic to $\mathbb{R}^{2n+1}$, and the multiplication
on it is defined by 
$$
 (x,y,z)(x',y',z')=\bigl(x+x',y+y',z+z'+\beta(x,y')\bigr),
$$
for $x,y,x',y'\in\mathbb{R}^n$ and $z,z'\in\mathbb{R}$.  Here $\beta(\ ,\ )$
is the usual inner product on $\mathbb{R}^n$.  We use this notation for a
possible future generalization.  For the extended Heisenberg Lie group
$\tilde{H}$ (corresponding to $\tilde{\Gh}$), see Example \ref{ex} in Appendix
below.
\end{defn}

 Taking advantage of the fact that their underlying spaces coincide, let us
from now on identify $H$ with $\Gh$ (as spaces) via the evident map: 
$$
(x,y,z) \mapsto {\sum_{i=1}^n\left(x_i\mathbf{x_i}+y_i\mathbf{y_i}\right)}
+z\mathbf{z},
$$
where $x=(x_1,\cdots,x_n)\in\mathbb{R}^n$, $y=(y_1,\cdots,y_n)\in\mathbb{R}^n$.
On this space $H\cong\Gh$, let us fix a Lebesgue measure.  This would be the
Haar measure for $H$.

 Note that this definition of the Heisenberg Lie group is different from the
one that is given by the Baker--Campbell--Hausdorff series for $\Gh$.  Thus
our set-up slightly differs (though isomorphic) from the setting in \S3 of
\cite{BJKp1}.  The use of the identification map as the diffeomorphism
between $\Gh$ and $H$ will make the subsequent calculation simpler.

 To obtain a Lie bialgebra structure on $\Gh$, consider $\Gg=\Gh^*$, the dual
vector space of $\Gh$, and fix a nonzero real number $\lambda$.  Let us define
the following Lie algebra structure on $\Gg$:
$$
[\mathbf{p_i},\mathbf{q_j}]=0,\quad 
[\mathbf{p_i},\mathbf{r}]=\lambda\mathbf{p_i},\quad [\mathbf{q_i},\mathbf{r}]
=\lambda\mathbf{q_i},
$$
where $\mathbf{p_i},\mathbf{q_i} (i=1,\dots,n),\mathbf{r}$ are the dual basis
for $\mathbf{x_i},\mathbf{y_i} (i=1,\dots,n),\mathbf{z}$.  Then we have the
following proposition:

\begin{prop}\label{bialgebra}
The (mutually dual) Lie algebras $\Gh$ and $\Gg$ determine a Lie bialgebra.
\end{prop}

\begin{proof}
We can prove this statement directly.  But let us choose an indirect method,
which would give us a deeper insight (and more information) about the situation.

 Consider the following element contained in $\tilde{\Gh}\otimes\tilde{\Gh}$:
\begin{equation}\label{(rmatrix)}
r=\lambda\bigl(\mathbf{z}\otimes\mathbf{d}+\mathbf{d}\otimes\mathbf{z}
+2\sum_{i=1}^n(\mathbf{x_i}\otimes\mathbf{y_i})\bigr).
\end{equation}
By elementary Lie algebra calculations, we can show that $r$ satisfies the
so-called ``classical Yang--Baxter equation'' (CYBE):
$$
[r_{12},r_{13}]+[r_{12},r_{23}]+[r_{13},r_{23}]=0.
$$
The notation $r_{ij}$ is understood as an element in $\tilde{\Gh}\otimes
\tilde{\Gh}\otimes\tilde{\Gh}$, and the meaning is fairly obvious (See
\cite{Dr}, \cite{CP}).  We can also show without difficulty that
$r_{12}+r_{21}$ is $\Gh$--invariant.  Therefore, $r$ is a ``(quasitriangular)
classical $r$--matrix'' (\cite{Dr}, \cite{CP}).

 Since we have a (quasitriangular) $r$--matrix, a ``coboundary'' Lie bialgebra
structure on $\tilde{\Gh}$ is defined by $\tilde{\delta}:\tilde{\Gh}\to
\tilde{\Gh}\wedge\tilde{\Gh}$, where
$$
\tilde{\delta}(X)=\operatorname{ad}_X(r).
$$
By restricting $\tilde{\delta}$ to $\Gh$, we obtain the map $\delta:\Gh\to
\Gh\wedge\Gh$, given by
$$
\delta(\mathbf{x_i})=\lambda\mathbf{x_i}\wedge\mathbf{z},\quad
\delta(\mathbf{y_i})=\lambda\mathbf{y_i}\wedge\mathbf{z},\quad
\delta(\mathbf{z})=0.
$$
This map is easily shown to be a 1--cocycle with respect to the adjoint
representation for $\Gh$, and hence $(\Gh,\delta)$ defines a Lie bialgebra.

 The Lie bialgebra structure $\delta$ on $\Gh$ determines a Lie bracket
on the dual vector space $\Gh^*$ by:
$$
\bigl\langle[\mu,\nu],X\bigr\rangle=\bigl\langle \mu\otimes\nu,\delta(X)
\bigr\rangle,
$$
where $\mu,\nu\in\Gh^*$, $X\in\Gh$, and $\langle\ ,\ \rangle$ is the dual
pairing between $\Gh^*$ and $\Gh$.  By straightforward calculation using
the definition of $\delta$, we can see that the resulting Lie bracket
coincides with the one we defined above on $\Gg=\Gh^*$.  This means that
the Lie bialgebra $(\Gh,\delta)$ is exactly the one determined by the pair
$\Gh$ and $\Gg$.
\end{proof}

\begin{cor}
By means of the classical $r$--matrix of \eqref{(rmatrix)} and the Lie
bialgebra $(\tilde{\Gh},\tilde{\delta})$ obtained from it, we can also
find the dual Lie bialgebra $\tilde{\Gg}=\tilde{\Gh}^*$ of $\tilde{\Gh}^*$:
It is spanned by the dual basis elements $\mathbf{p_i},\mathbf{q_i} (i=1,
\dots,n),\mathbf{r},\mathbf{s}$, with $\mathbf{p_i},\mathbf{q_i},\mathbf{r}$
satisfying the same relation as before and $\mathbf{s}$ being central.  By
construction, $(\Gh,\delta)$ or $(\Gh,\Gg)$ is a ``sub-bialgebra'' of $(\tilde
{\Gh},\tilde{\delta})$ or $(\tilde{\Gh},\tilde{\Gg})$.
\end{cor}

\begin{rem}
Unlike $(\tilde{\Gh},\tilde{\delta})$, the Lie bialgebra $(\Gh,\delta)$ cannot
be obtained as a coboundary from any classical $r$--matrix contained in $\Gh
\otimes\Gh$.  Thus the introduction of the extended Heisenberg Lie algebra
$\tilde{\Gh}$ is essential.  The same situation occurs in \cite{ALT}, \cite{CGST},
where the authors find (using the same classical $r$--matrix as above) a
quantized universal enveloping algebra (i.\,e. QUE algebra) deformation of
the Heisenberg Lie algebra.
\end{rem}

 The Lie group $G$ associated with $\Gg$ is, by definition, the {\em dual 
Poisson--Lie group\/} of $H$.  To know more about $G$, note first that the
Lie algebra $\Gg$ is a semi-direct product of its two (abelian) subalgebras
$\Gm=\operatorname{span}(\mathbf{r})$ and $\Gq=\operatorname{span}(\mathbf{p_i},
\mathbf{q_i}|i=1,\dots,n)$.  This is evident from its defining relations.
Therefore, the connected and simply connected Lie group $G$ associated with
$\Gg$ should be a semi-direct product group.  Since $\Gm$ and $\Gq$ are abelian
Lie algebras, they are identified (as spaces) with their corresponding abelian
Lie groups.  This suggests the following definition of $G$:

\begin{defn}\label{dual group}(The dual Poisson--Lie group)
Let $G=\Gq\times\Gm=\Gg$ as a vector space. Define the multiplication law
on it by
$$
(p,q,r)(p',q',r')=(e^{\lambda r'}p+p',e^{\lambda r'}q+q',r+r').
$$
\end{defn}

 Here, $p=(p_1,\dots,p_n)\in\mathbb{R}^n$, $q=(q_1,\dots,q_n)\in\mathbb{R}^n$,
$r\in\mathbb{R}$, and we are identifying $(p,q,r)\in G$ with the element
${\sum_{i=1}^n(p_i\mathbf{p_i}+q_i\mathbf{q_i})}+r\mathbf{r}\in\Gg$.  This
means that, in particular, $G$ is an exponential solvable Lie group.  Now
on $G$ (which is being identified as a space with $\Gg=\Gh^*$), let us choose
the Plancherel Lebesgue measure dual to the fixed Haar measure on $H\cong\Gh$
(See the remark made below, at the end of this section.).  This will be the
{\em left invariant\/} Haar measure on $G$.

 The group $G$ will be our main object of study:  Following the method of
\cite{BJKp1}, we are going to find a deformation quantization of $G$, using
the duality between $\Gh$ and $\Gg$ (or $H$ and $G$).  Before we begin our
main discussion, let us make a short remark on Fourier transforms between
dual spaces.  This will serve a purpose of setting up the notation we will
be using in this paper.

\begin{rem}(Fourier transforms between dual spaces, Plancherel measure)
 Let $W$ be a (real) vector space.  Let us fix a Lebesgue measure, $dx$, on
$W$.  Let $W^*$ be the dual vector space of $W$.  We choose on $W^*$ the dual
``Plancherel measure'', denoted by $d\mu$, which is also a Lebesgue  measure. 
Then the Fourier transform from $L^2(W)$ to $L^2(W^*)$ is given by
$$
({\mathcal F}\xi)(\mu)
 =\int_W\bar{e}\bigl[\langle\mu,x\rangle\bigr]\xi(x)\,dx.
$$
Here $\langle\ ,\ \rangle$ denotes the dual pairing between $W^*$ and $W$,
and $e(\ )$ is the function defined by $e(t)=e^{2\pi it}$.  So $\bar{e}(t)
=e^{-2\pi it}$.  By our choice of measures, the Fourier transform is a 
unitary operator, whose inverse is the following inverse Fourier transform:
$$
({\mathcal F}^{-1}\zeta)(x)
 =\int_{W^*}e\bigl[\langle\mu,x\rangle\bigr]\zeta(\mu)\,d\mu.
$$
If $Z$ is a subspace of $W$ and if we fix a Lebesgue measure, $dz$, on $Z$,
there is a unique Plancherel Lebesgue measure, $d\dot{x}$, on $W/Z$ so that
$dx=d\dot{x}dz$.  Since $Z^{\bot}=(W/Z)^*$, we can also choose as above an
appropriate Plancherel measure, $dq$, on $Z^{\bot}\subseteq W^*$.  This
enables us to define the ``partial'' Fourier transform from $L^2(W/Z
\times W^*/{Z^{\bot}})$ to $L^2(W^*)=L^2(Z^{\bot}\times W^*/{Z^{\bot}})$,
given by
$$
f^{\wedge}(q,r)=\int_{W/Z}\bar{e}\bigl[\langle q,\dot{x}\rangle\bigr]
f(\dot{x},r)\,d\dot{x}.
$$
Its inverse, $\phi\mapsto\phi^{\vee}$, is defined similarly as above.

 Let $S(W)$ denote the space of Schwartz functions on $W$.  Then by Fourier
transform, $S(W)$ is carried onto $S(W^*)$ and vice versa.  The Fourier
inversion theorem (the unitarity of the Fourier transform) implies that
we have: ${\mathcal F}^{-1}({\mathcal F}f)=f$ for $f\in S(W)$ and  
${\mathcal F}({\mathcal F}^{-1}\phi)=\phi$ for $\phi\in S(W^*)$.  Similar
assertion is true for the partial Fourier transform.
\end{rem}

\section{Deformation quantization of $G$}

 Let us compute explicitly the compatible Poisson bracket on the dual Poisson--Lie
group $G$.  Later in this section, we are going to find a deformation quantization
of $G$ in the direction of this Poisson bracket.  To compute the Poisson bracket,
let us first compute the Lie bialgebra structure $(\Gg,\theta)$ on $\Gg$.  Since
$\theta$ determines the dual Lie bialgebra of $(\Gh,\delta)$, it should be the
dual map of the given Lie bracket on $\Gg^*=\Gh$:

\begin{lem}
Let $\theta:\Gg\to\Gg\wedge\Gg$ be defined by its values on the basis
vectors of $\Gg$ as follows:
$$
\theta(\mathbf{p_i})=0,\quad\theta(\mathbf{q_i})=0,\quad
\theta(\mathbf{r})=\sum_{i=1}^n(\mathbf{p_i}\otimes\mathbf{q_i}
-\mathbf{q_i}\otimes\mathbf{p_i})=\sum_{i=1}^n(\mathbf{p_i}
\wedge\mathbf{q_i}).  
$$
Then $\theta$ is the dual map of the Lie bracket on $\Gh$.  Hence it is the 
1--cocycle giving the dual Lie bialgebra structure on $\Gg$.
\end{lem}

\begin{proof}
Straightforward.
\end{proof}

 By using the simple connectedness of $G$, the Lie bialgebra structure
$(\Gg,\theta)$ determines the compatible Poisson bracket on $G$  (See
\cite{LW}, \cite{CP}.).  The calculation and the result is given below.
See \cite{SZ} for a similar result.  Observe also that our expression of
the Poisson bracket does not depend on the $p$ and $q$ variables.  

\begin{theorem}\label{poissonbr}
The Poisson bracket on the dual Poisson--Lie group $G$ is given by
the following expression:  For $\phi,\psi\in C^{\infty}(G)$,  
\begin{equation}\label{(2.1)}
\{\phi,\psi\}(p,q,r)=\left(\frac{e^{2\lambda r}-1}{2\lambda}\right)
\bigl(\beta(x,y')-\beta(x',y)\bigr), \qquad (p,q,r)\in G
\end{equation}
where $d\phi(p,q,r)=(x,y,z)$ and $d\psi(p,q,r)=(x',y',z')$, which are
naturally considered as elements of $\Gh$.
\end{theorem}

\begin{proof}
Let $\operatorname{Ad}:G\to\operatorname{Aut}(\Gg)$ be the adjoint 
representation of $G$ on $\Gg$.  We have to look for a map $F:G\to\Gg
\wedge\Gg$, which is a group 1--cocycle on $G$ for $\operatorname{Ad}$ 
and whose derivative at the identity element, $dF_e$, coincides with 
the map $\theta$.  Since $\theta$ depends only on the $r$--variable, so
should $F$.  Thus we only need to look for a map $F$ satisfying the 
condition:
$$
F(r_1+r_2)=F(r_1)+\operatorname{Ad}_{(0,0,r_1)}\bigl(F(r_2)\bigr),
$$
such that its derivative at the identity element is the map, 
$dF_e(r)=\theta(r)=r{\sum_{i=1}^n}(\mathbf{p_i}\wedge\mathbf{q_i})$. 
Meanwhile, note that the representation $\operatorname{Ad}$ sends the 
basis vectors of $\Gg$ as follows:
\begin{align}
\operatorname{Ad}_{(0,0,r')}(\mathbf{p_i})&=(0,0,r')(1,0,0)(0,0,-r')
=(e^{-\lambda r'},0,0)=e^{-\lambda r'}\mathbf{p_i}, \notag \\
\operatorname{Ad}_{(0,0,r')}(\mathbf{q_i})&=e^{-\lambda r'}\mathbf{q_i},
\quad \operatorname{Ad}_{(0,0,r')}(\mathbf{r})=\mathbf{r}. \notag
\end{align}
So the 1--cocycle condition for $F$ becomes:
$$
F(r_1+r_2)=F(r_1)+e^{-2\lambda r_1}F(r_2).
$$
From this equation together with the condition, $dF_e=\theta$, we 
obtain:
$$
F(p,q,r)=F(r)=\left(\frac{1-e^{-2\lambda r}}{2\lambda}\right)
{\sum_{i=1}^n}(\mathbf{p_i}\wedge\mathbf{q_i}).
$$

The Poisson bivector field is the right translation of this 1--cocycle
$F$, given by ${R_{(p,q,r)}}_*F(p,q,r)$.  Since the right translations
are ${R_{(p,q,r)}}_*(\mathbf{p_i})=e^{\lambda r}\mathbf{p_i}$ and 
${R_{(p,q,r)}}_*(\mathbf{q_i})=e^{\lambda r}\mathbf{q_i}$, we obtain 
equation \eqref{(2.1)} for our Poisson bracket by:
$$
\{\phi,\psi\}(p,q,r)=\bigl\langle{R_{(p,q,r)}}_*F(p,q,r),d\phi(p,q,r)
\wedge d\psi(p,q,r)\bigr\rangle.
$$
\end{proof}

 Since we will use the expression $(e^{2\lambda r}-1)/{2\lambda}$ quite
often, let us give it a special notation, $\eta_{\lambda(r)}$.  This function
satisfies a convenient identity, which is given in Lemma \ref{lem}.  The proof
is straightforward.

\begin{defn}
Let $\lambda\in\mathbb{R}$ be fixed.  Let us denote by $\eta_{\lambda}$
the function on $\mathbb{R}$ defined by
$$
\eta_{\lambda}(r)=\frac{e^{2\lambda r}-1}{2\lambda}.
$$
When $\lambda=0$, we define $\eta_0(r)=r$.
\end{defn}

\begin{lem}\label{lem}
For $r,r'\in\Gg/\Gq$, we have:
\begin{equation}\label{(lemeq)}
(e^{-2\lambda r'})\eta_{\lambda}(r+r')-(e^{-2\lambda r'})\eta_{\lambda}
(r')=\eta_{\lambda}(r).
\end{equation}
\end{lem}

 Since we are identifying $G\cong\Gg$ as spaces, our Poisson bracket on $G$
may as well be regarded as a Poisson bracket on $\Gg=\Gh^*$.  It is a non-linear
Poisson bracket, but it is nevertheless of the special type studied in
\cite{BJKp1}.  We summarize this observation in the next proposition.  Here
$\Gz$ denotes the center of $\Gh$, spanned by $\mathbf{z}$.  Also $\Gq=
\Gz^{\bot}$, in $\Gg$.  As before, we regard the vectors $x,y,x',y'\in
\mathbb{R}^n$ as elements of $\Gh/\Gz=\operatorname{span}(\mathbf{x_i},
\mathbf{y_i}|i=1,\dots,n)$, and similarly $r\in\mathbb{R}$ as an element of
$\Gg/\Gq$.

\begin{prop}
\begin{enumerate}
\item Let $\omega:\Gh/\Gz\times\Gh/\Gz\to C^{\infty}(\Gg/\Gq)$ be the map
defined by
$$
\omega\bigl((x,y),(x',y');r\bigr)=\eta_{\lambda}(r)\bigl(\beta(x,y')
-\beta(x',y)\bigr).
$$
Then it is a Lie algebra cocycle for $\Gh/\Gz$ having values in 
$V=C^{\infty}(\Gg/\Gq)$, regarded as a trivial $U(\Gh/\Gz)$-module.
\item The Poisson bracket on $\Gg=\Gh^*$ given by equation \eqref{(2.1)} is
realized as a sum of the (trivial) linear Poisson bracket on $(\Gh/\Gz)^*$
and $\omega$.  
\item The space $V=C^{\infty}(\Gg/\Gq)$ is canonically contained in 
$C^{\infty}(\Gg)$ such that $\Gh\cap V=\Gz$.
\end{enumerate}
Thus we conclude that our Poisson bracket is the ``cocycle perturbation''
(in the sense of \cite{BJKp1}) of the linear Poisson bracket on $\Gh^*$.
\end{prop}

\begin{proof}
We can see easily that $\omega$ is a skew-symmetric, bilinear map, trivially
satisfying the cocycle identity since $\Gh/\Gz$ is abelian.  Since $\Gh/\Gz$
is an abelian Lie algebra, it also follows that the linear Poisson bracket
on $(\Gh/\Gz)^*$ is the trivial one.  Thus the second assertion of the 
proposition is immediate from the definition of $\omega$.

The functions in $V=C^{\infty}(\Gg/\Gq)$ can be canonically realized as 
functions in $C^{\infty}(\Gg)$ by the ``pull-back'' using the natural
projection of $\Gg$ onto $\Gg/\Gq$.  If we regard the elements of $\Gh$ also
as (linear) functions in $C^{\infty}(\Gg)$, we have: $\Gh\cap V=\Gz$.  It 
follows that our Poisson bracket is an extension of the linear Poisson 
bracket on $(\Gh/\Gz)^*$ by the cocycle $\omega$.  We showed in \cite{BJKp1} 
(See Theorems 2.2 and 2.3) that this formulation is equivalent to viewing
the Poisson bracket as a ``cocycle perturbation'' of the linear Poisson
bracket on $\Gh^*$.
\end{proof}

\begin{rem}
When $\lambda=0$, we have:
$$
\omega_0\bigl((x,y),(x',y');r\bigr)=r\bigl(\beta(x,y')-\beta(x',y)\bigr).
$$
It is a linear function on $\Gg/\Gq$, so we may write it as:
\begin{equation}\label{(omega0)}
\omega_0\bigl((x,y),(x',y')\bigr)=\bigl(\beta(x,y')-\beta(x',y)\bigr)
\mathbf{z}.
\end{equation}
Thus $\omega_0$ is a cocycle for $\Gh/\Gz$ having values in $\Gz$.  It is
clear that the linear Poisson bracket on $\Gg=\Gh^*$ (see \cite{SZ}) is
determined by the cocycle $\omega_0$.  In other words, the ``perturbation''
is given by (nonzero) $\lambda$ and the associated cocyle $\omega$.
\end{rem}

 Deformation quantization of our Poisson bracket on $\Gg$, which we will denote
by $\{\ ,\ \}_{\omega}$ from now on, is obtained by following the steps of
\cite{BJKp1}.  First, we construct from the given Lie algebra cocycle $\omega$
the continuous family of $\mathbb{T}$--valued group cocycles for the Lie group
$H/Z$ of $\Gh/\Gz$.  Here $Z=\operatorname{span}\{\mathbf{z}\}$ is the Lie
subgroup of $H$ corresponding to $\Gz$.

\begin{prop}
Consider the map $R:H/Z\times H/Z\to V=C^{\infty}(\Gg/\Gq)$ defined by
$$
R\bigl((x,y),(x',y');r\bigr)=\eta_{\lambda}(r)\beta(x,y').
$$
Then $R$ is a group cocycle for $H/Z$ having values in $V$, regarded as
an additive abelian group.  Fix now an element $r\in\Gg/\Gq$.  Define 
the map $\sigma^{r}:H/Z\times H/Z\to\mathbb{T}$ by
$$
\sigma^{r}\bigl((x,y),(x',y')\bigr)=\bar{e}\bigl[R((x,y),(x',y');r)
\bigr]=\bar{e}\bigl[\eta_{\lambda}(r)\beta(x,y')\bigr].
$$
Then each $\sigma^{r}$ is a $\mathbb{T}$--valued, normalized group
cocycle for $H/Z$.  And $r\mapsto\sigma^{r}$ forms a continuous field
of cocycles.
\end{prop}

\begin{proof}
Let $h=(x,y),h'=(x',y'),h''=(x'',y'')$ be elements of $H/Z$.  
Then for $r\in\Gg/\Gq$, we have the cocycle identity:
\begin{align}
\sigma^{r}(hh',h'')\sigma^{r}(h,h')&=\bar{e}\bigl[\eta_{\lambda}(r)
(\beta(x,y'')+\beta(x',y'')+\beta(x,y'))\bigr] \notag \\
&=\sigma^{r}(h,h'h'')\sigma^{r}(h',h'').  \notag  
\end{align}
We also have: $\sigma^{r}(h,0)=1=\sigma^{r}(0,h)$, where $0=(0,0)$ is
the identity element of $H/Z$.  From the definition, the continuity is
also clear.
\end{proof}

 Let us consider $S_{3c}(\Gh/\Gz\times\Gg/\Gq)$, the space of Schwartz functions
in the $(x,y,r)$--variables having compact support in the $r\in\Gg/\Gq$ variable.
Since $\Gh/\Gz$ is identified with the abelian group $H/Z$, we can regard $S_{3c}
(\Gh/\Gz\times\Gg/\Gq)$ as contained in $L^1\bigl(H/Z,C_{\infty}(\Gg/\Gq)\bigr)$.
Using the continuous field of cocycles $\sigma$, we can define on it the following
twisted convolution:
$$
(f*_{\sigma}g)(x,y,r)=\int f(\tilde{x},\tilde{y},r)g(x-\tilde{x},
y-\tilde{y},r)\bar{e}\bigl[\eta_{\lambda}(r)\beta(\tilde{x},y-\tilde{y})
\bigr]\,d\tilde{x}d\tilde{y}.
$$
It is not difficult to see that $S_{3c}(\Gh/\Gz\times\Gg/\Gq)$ is indeed an
algebra.

 To transfer this algebra structure to the level of functions on $\Gg$, we
introduce the partial Fourier transform.  The partial Fourier transform,
${}^{\wedge}$, from $S(\Gh/\Gz\times\Gg/\Gq)$ to $S(\Gg)=S(\Gq\times\Gg/\Gq)$
is defined by  
$$
f^{\wedge}(p,q,r)=\int\bar{e}(p\cdot x+q\cdot y)f(x,y,r)\,dxdy,
$$
where $p\cdot x+q\cdot y$ is the dual pairing between $(p,q)\in\Gq$ and
$(x,y)\in\Gh/\Gz$.  The inverse partial Fourier transform, ${}^{\vee}$, from
$S(\Gg)$ to $S(\Gh/\Gz\times\Gg/\Gq)$ is defined in a similar manner, with
$\bar{e}(\ )$ replaced by $e(\ )$.  We are again assuming that we have chosen
appropriate Plancherel measures for $\Gh/\Gz$ and $\Gq=\Gz^{\bot}$, so that
the Fourier inversion theorem is valid.

 To define the deformed multiplication between the functions on $\Gg$, consider
the subspace ${\mathcal A}=S_{3c}(\Gg)\subseteq S(\Gg)$, which is the image
under the partial Fourier transform, ${}^{\wedge}$, of the twisted convolution
algebra $S_{3c}(\Gh/\Gz\times\Gg/\Gq)$.

\begin{prop}\label{product}
Let ${\mathcal A}=S_{3c}(\Gg)$ be the space of Schwartz functions on $\Gg$
having compact support in the $r$--variable.  On $\mathcal A$, define the
deformed multiplication, $\times$, by: $\phi\times\psi=(\phi^{\vee}\ast
_{\sigma}\psi^{\vee})^{\wedge}$, for $\phi,\psi\in{\mathcal A}$.  We then
obtain:
$$
(\phi\times\psi)(p,q,r)=\int\bar{e}\bigl[(p-p')\cdot\tilde{x}\bigr]
\phi(p',q,r)\psi(p,q+\eta_{\lambda}(r)\tilde{x},r)\,dp'd\tilde{x}. 
$$
\end{prop}

\begin{proof}
Use the Fourier inversion theorem to the expression:
\begin{align}
&(\phi\times\psi)(p,q,r)=(\phi^{\vee}\ast_{\sigma}\psi^{\vee})^{\wedge}
(p,q,r)   \notag\\
&\quad=\int\bar{e}(p\cdot x+q\cdot y)\phi^{\vee}(\tilde{x},\tilde{y},r)
\psi^{\vee}(x-\tilde{x},y-\tilde{y},r)\bar{e}\bigl[\eta_{\lambda}(r)
\beta(\tilde{x},y-\tilde{y})\bigr]  \notag \\
&\quad\quad\quad\quad d\tilde{x}d\tilde{y}dxdy \notag
\end{align}
\end{proof}

 Note that when $\lambda=0$, the operation $\ast_{\sigma}$ on $L^1\bigl(H/Z,
C_{\infty}(\Gg/\Gq)\bigr)$ is given by the cocycle, $\bigl((x,y),(x',y')\bigr)
\mapsto\bar{e}\bigl[r\beta(x,y')\bigr]$ for $H/Z$.  But this is essentially
the ordinary convolution on $S(H)$ transferred to the functions in the
$(x,y,r)$--variables.  Compare this with our case, with the cocycle, ${\sigma}^r
\bigl((x,y),(x',y')\bigr)=\bar{e}\bigl[\eta_{\lambda}(r)\beta(x,y')\bigr]$.
We can see that the passage from the linear Poisson bracket (when $\lambda=0$)
to our ``perturbed'' (non-linear) Poisson bracket corresponds to the ``change
of cocycles'', or to the passage from ordinary convolution to twisted convolution.

 The situations between linear Poisson bracket case (\cite{Rf3}) and our
perturbed case are quite similar, and this will be exploited from time
to time.  However, in our more general case, the space $S(\Gg)$ is no
longer closed under the deformed multiplication.  We had to define the
multiplication in its subspace ${\mathcal A}$.

 The algebra ${\mathcal A}$ is shown to be a pre-$C^*$--algebra, whose 
involution and the $C^*$--norm are again obtained using the partial Fourier
transform between ${\mathcal A}$ and $S_{3c}(\Gh/\Gz\times\Gg/\Gq)$, the 
latter being viewed as a (dense) subalgebra of the ${}^*$--algebra $L^1\bigl
(H/Z,C_{\infty}(\Gg/\Gq),\sigma\bigr)$.  See \cite{BJKp1}, for the precise
definitions of the ${}^*$--algebra operations.

\begin{prop}\label{invol}
Let $\mathcal A$ be as above.
\begin{enumerate}
\item The involution on $\mathcal A$ is defined by: $\phi\mapsto\bigl
((\phi^{\vee})^*\bigr)^{\wedge}$, where ${}^*$ denotes the involution on
$S_{3c}(\Gh/\Gz\times\Gg/\Gq)$.  If we denote the involution on $\mathcal A$
by the same notation, ${}^*$, then we have:
$$
\phi^*(p,q,r)=\int\overline{\phi(p',q',r)}\bar{e}\bigl[(p-p')\cdot x+
(q-q')\cdot y\bigr]\bar{e}\bigl[\eta_{\lambda}(r)\beta(x,y)\bigr]\,dp'dq'dxdy
$$
\item Via partial Fourier transform, we also define the canonical $C^*$--norm
on $\mathcal A$, by transferring the canonical $C^*$--norm on $S_{3c}(\Gh/\Gz
\times\Gg/\Gq)\subseteq L^1\bigl(H/Z,C_{\infty}(\Gg/\Gq),\sigma\bigr)$.
\end{enumerate}
\end{prop}

\begin{proof}
The involution on $S_{3c}(\Gh/\Gz\times\Gg/\Gq)$ is given by
$$
f^*(x,y,r)=\overline{f(-x,-y,r)}\bar{e}\bigl[\eta_{\lambda}(r)\beta(x,y)\bigr].
$$
It is easy to see that $S_{3c}(\Gh/\Gz\times\Gg/\Gq)$ is closed under
the involution.  We transfer this operation to the $\mathcal A$ level by
$\phi\mapsto\bigl((\phi^{\vee})^*\bigr)^{\wedge}$.  Use the Fourier
inversion theorem to obtain the above expression.

 On $L^1\bigl(H/Z,C_{\infty}(\Gg/\Gq),\sigma\bigr)$, one has the canonical
$C^*$--norm, $\|\ \|$, such that the completion with respect to $\|\ \|$ of
this $L^1$ algebra is the enveloping $C^*$--algebra $C^*\bigl(H/Z,C_{\infty}
(\Gg/\Gq),{\sigma}\bigr)$, called the twisted group $C^*$--algebra.  By
$\phi\mapsto\|\phi^{\vee}\|$, we can define on $\mathcal A$ its $C^*$--norm,
still denoted by $\|\ \|$.
\end{proof}

 Since the function space $S_{3c}(\Gh/\Gz\times\Gg/\Gq)$ is dense in $L^2
(\Gh/\Gz\times\Gg/\Gq)$ with respect to the $L^2$--norm, its product
$\ast_{\sigma}$ corresponds to a representation of $S_{3c}(\Gh/\Gz\times
\Gg/\Gq)$ on $L^2(\Gh/\Gz\times\Gg/\Gq)$ such that the functions acts as
the multiplication operators.  This representation is naturally extended
to $C^*(H/Z,C_{\infty}(\Gg/\Gq),\sigma)$.  More precisely, we have a
representation, $L$, of the twisted group $C^*$--algebra on $L^2(\Gh/\Gz
\times\Gg/\Gq)$ defined by
$$
(L_f\xi)(x,y,r)=\int f(x',y',r)\xi(x-x',y-y',r)\bar{e}\bigl[\eta_
{\lambda}(r)\beta(x',y-y')\bigr]\,dx'dy',
$$
for $f\in S_{3c}(\Gh/\Gz\times\Gg/\Gq)$ and $\xi\in L^2(\Gh/\Gz\times
\Gg/\Gq)$.  It is actually a {\em (left) regular representation\/} of 
$C^*(H/Z,C_{\infty}(\Gg/\Gq),\sigma)$, induced from a (faithful) 
representation of $C_{\infty}(\Gg/\Gq)$ on $L^2(\Gg/\Gq)$ given by 
multiplication.

 In what follows, we will be working with the Hilbert space $L^2(\Gh/\Gz
\times\Gg/\Gq)$ most of the time.  So let us from now on denote this
Hilbert space by $\mathcal H$.  Via the isomorphism between $S_{3c}(\Gh/\Gz
\times\Gg/\Gq)$ and $\mathcal A$, the representation $L$ may as well be
regarded as a representation of $\mathcal A$ on $\mathcal H$.  Let us also
denote this representation by $L$.  Then for $\phi\in{\mathcal A}$, 
\begin{align}
&(L_{\phi}\xi)(x,y,r)  \notag \\
&=\int\phi^{\vee}(x',y',r)\xi(x-x',y-y',r)
\bar{e}\bigl[\eta_{\lambda}(r)\beta(x',y-y')\bigr]\,dx'dy'  \notag \\
&=\int\phi(p,q,r)\xi(x-x',y-y',r)e(p\cdot x'+q\cdot y') 
\bar{e}\bigl[\eta_{\lambda}(r)\beta(x',y-y')\bigr]\,dpdqdx'dy'. \label{(rep)}
\end{align}
It is clear that $L$ is equivalent to the representation of $\mathcal A$ on
$L^2(\Gg)$ given by the multiplication $\times$.  The partial Fourier transform
is the intertwining unitary operator between the Hilbert spaces $\mathcal H$
and $L^2(\Gg)$. 

 The representation $L$ is the regular representation induced from a faithful
representation of $C^{\infty}(\Gg/\Gq)$.  So the corresponding $C^*$--norm 
and the completion will give us the ``reduced'' twisted group $C^*$--algebra 
$C^*_r\bigl(H/Z,C_{\infty}(\Gg/\Gq),{\sigma}\bigr)$.  Since $H/Z$ is abelian,
the amenability condition holds in our case.  i.\,e. $C^*_r\bigl(H/Z,
C_{\infty}(\Gg/\Gq),{\sigma}\bigr)\cong C^*\bigl(H/Z,C_{\infty}(\Gg/\Gq),
{\sigma}\bigr)$.  This follows from the result of Packer and Raeburn \cite{PR1},
which says that the amenability of the group implies the amenability of the
twisted group $C^*$--algebra.  Because of the amenability, we can see that
for $\phi\in{\mathcal A}$, we have: $\|\phi\|=\|L_{\phi}\|$.

\begin{defn}\label{cstar}
Let $\mathcal A$ be defined as above and let it be equipped with the multiplication,
$\times$, given by Proposition \ref{product} and the involution, ${}^*$, given by
Proposition \ref{invol}.  Let us denote by $A$ the $C^*$--completion of $\mathcal A$
with respect to the norm defined by $\|\phi\|=\|L_{\phi}\|$, where $L_{\phi}$ is
regarded as an operator on $\mathcal H$ by equation \eqref{(rep)}.  This is the
$C^*$--algebra we will be interested in throughout the rest of this paper.  It is
clear that $A\cong C^*_r\bigl(H/Z,C_{\infty}(\Gg/\Gq),{\sigma}\bigr)\cong 
C^*\bigl(H/Z,C_{\infty}(\Gg/\Gq),{\sigma}\bigr)$.
\end{defn}

 Recall that we are identifying $G$ with $\Gg$ as spaces and the Plancherel
Lebesgue measure on $\Gg$ we have been using coincides with the Haar measure
on $G$.  We thus have, as a (dense) subspace, ${\mathcal A}\subseteq C_{\infty}
(G)$.  And the results we obtained so far about functions on $\Gg$ hold true
for functions on $G$.  The deformed function algebra $({\mathcal A},\times,{}^*)$,
as well as its $C^*$--completion $A$, provide a deformation quantization (to
be made precise shortly) of $C_{\infty}(G)$.

 At each of the steps above, we could have kept the parameter $\hbar$ as in
\cite{BJKp1}.  In our case, the deformed algebra would be isomorphic to the
twisted group algebra of $(H/Z)_{\hbar}=H/Z$ with the cocycle $\sigma_{\hbar}$
given by ${\sigma}_{\hbar}^r\bigl((x,y),(x',y')\bigr)=\bar{e}\bigl[\hbar
\eta_{\lambda}(r)\beta(x,y')\bigr]$.  Since $H/Z$ is abelian, the group doesn't
have to vary and only the cocycle ${\sigma}$ varies under the introduction of
the parameter $\hbar$.  See \cite{BJKp1} for more precise formulation.

 Let us denote by $(\times_{\hbar},{}^{*_{\hbar}},\|\ \|_{\hbar})$ the
corresponding operations on ${\mathcal A}$ obtained by the introduction
of the parameter $\hbar$.  The above discussion means that all we have to do
is replace $\beta(\ ,\ )$ by $\hbar\beta(\ ,\ )$.  Then define $A_{\hbar}$
as the $C^*$--completion of ${\mathcal A}$ with respect to $\|\ \|_{\hbar}$.
By the main theorem (Theorem 3.4) of \cite{BJKp1}, we thus obtain a (strict)
deformation quantization of our Poisson bracket $\{\ ,\ \}_{\omega}$ on $G$.

\begin{theorem}\label{main}
Consider the dual Poisson bracket on $G\cong\Gg$ defined by equation \eqref
{(2.1)}.  Let ${\mathcal A}=S_{3c}(G)$ be the subspace of $S(G)$ defined above.
For any $\hbar\in\mathbb{R}$, define a deformed multiplication and an involution
on $\mathcal A$, and also a $C^*$--norm on it, by replacing $\beta(\ ,\ )$ with
$\hbar\beta(\ ,\ )$ in Definition \ref{cstar}. Then $\bigl({\mathcal A},
\times_{\hbar},{}^{*_{\hbar}},\|\ \|_{\hbar}\bigr)_{\hbar\in\mathbb{R}}$
provides a strict deformation quantization (in the sense of \cite{Rf1, Rf4}) of
$\mathcal A$ in the direction of $(1/{2\pi})\{\ ,\ \}_{\omega}$.  In particular,
we have:
\begin{equation}\label{(deform)}
\lim_{\hbar\to 0}\left\|\frac{\phi\times_{\hbar}\psi-\psi\times_{\hbar}
\phi}{\hbar}-\frac i{2\pi}\{\phi,\psi\}_{\omega}\right\|_\hbar=0. 
\end{equation}
\end{theorem}

\begin{proof}
For full proof of the theorem, refer to Theorem 3.4 in \cite{BJKp1}, of which
ours is a special case.  We will briefly mention here a few of the main points
of proof.

First, we have to show that the family of $C^*$--algebras $\{A_{\hbar}\}
_{\hbar\in\mathbb{R}}$, where each $A_{\hbar}$ is the $C^*$--completion
of $\mathcal A$ with respect to $\|\ \|_{\hbar}$, forms a continuous
field of $C^*$--algebras.  Since each $A_{\hbar}$ is essentially a 
twisted group $C^*$--algebra of an abelian group $H/Z$, and only the
cocycle is being changed, the proof is actually simpler than in
\cite{BJKp1}.

Second, to prove the deformation property, it suffices to show that on
${\mathcal A}$, the expression, $(\phi\times_{\hbar}\psi-\psi\times_{\hbar}
\phi)/{\hbar}-(i/{2\pi})\{\phi,\psi\}_{\omega}$ has an $L^1$--bound.  Then
by Lebesgue's dominated convergence theorem, we would have the convergence in
the $L^1$--norm, which in turn gives  the convergence \eqref{(deform)} since
the $L^1$--norm dominates all the $C^*$--norms $\|\ \|_{\hbar}$.  The proof
crucially uses the fact that our functions are Schwartz functions having
compact support in the $r\in\Gg/\Gq$ variable.
\end{proof}

 From now on, we will fix the parameter $\hbar$ (e.\,g. $\hbar=1$) and 
take the resulting $C^*$--algebra $A$ as the candidate for our quantum group.  
If we want to specify the deformation process, we can always re-introduce
$\hbar$, and follow the arguments above.

 To summarize, the meaning of above construction and Definition \ref{cstar}
is that we are viewing the functions in ${\mathcal A}\subseteq S(G)$ as
operators on $\mathcal H$, by the regular representation $L$.  This naturally
defines the deformed multiplcation on $\mathcal A$, which is shown to be a
deformation quantization of $\bigl(G,\{\ \}_{\omega}\bigr)$ by Theorem \ref{main}.
So from now on, we will interprete $\phi$ and $L_{\phi}$ as the same.  Viewing
$\phi\in{\mathcal A}$ as a function has an advantage when we try to establish
a correspondence between our quantum setting and the classical, Poisson--Lie
group level.  While, viewing it as an operator $L_{\phi}\in A\subseteq
{\mathcal B}({\mathcal H})$ is essential to make our discussion rigorous
at the $C^*$--algebra level of ``locally compact quantum groups''.

 Meanwhile, observe that $L_{\phi}$ can be written as
\begin{equation}\label{(lphi)}
L_{\phi}=\int_H({\mathcal F}^{-1}\phi)(x,y,z)L_{x,y,z}\,dxdydz,
\end{equation}
where ${\mathcal F}^{-1}$ is the (inverse) Fourier transform from $\mathcal A$
into $S(\Gh)$, and $L_{a,b,c}$ for $(a,b,c)\in H$ is the operator on $\mathcal H$
defined by
\begin{equation}\label{(block)}
(L_{a,b,c}\xi)(x,y,r)=\bar{e}(rc)\bar{e}\bigl[\eta_{\lambda}(r)\beta(a,y-b)\bigr]
\xi(x-a,y-b,r).
\end{equation}

 By using the Fourier inversion theorem purely formally to this expression,
$L_{a,b,c}$ can be written as:
\begin{align}
(L_{a,b,c}\xi)(x,y,r)&=\int\bar{e}\bigl[\langle(p,q,r),(a,b,c)\rangle\bigr]
e(p\cdot\tilde{x}+q\cdot\tilde{y})  \notag \\ 
&\qquad\quad \bar{e}\bigl[\eta_{\lambda}(r)\beta(\tilde{x},y-\tilde{y})
\bigr]\xi(x-\tilde{x},y-\tilde{y},r)\,dpdqd\tilde{x}d\tilde{y}. \notag
\end{align}
Comparing this with equation \eqref{(rep)}, we observe that $L_{a,b,c}$ can be
regarded as a (continuous) function on $G$ defined by
$$
L_{a,b,c}(p,q,r)=\bar{e}\bigl[\langle(p,q,r),(a,b,c)\rangle\bigr]
=\bar{e}[p\cdot a+q\cdot b+rc].
$$
Note however that $L_{a,b,c}$ is not contained in $\mathcal A$.  It is not even
an element of $A$.  A more precise statement is that $L_{a,b,c}$ is a multiplier
(i.\,e. an element of $M(A)$).

 By \eqref{(lphi)}, any representation $\Pi$ of ${\mathcal A}$ or $A$ will be
written as
$$
\Pi(\phi)=\Pi(L_{\phi})=\int_H({\mathcal F}^{-1}\phi)(x,y,z)\Pi(L_{x,y,z})\,dxdydz.
$$
This means that if we have to check whether two non-degenerate representations
are equal, it suffices to check that they are equal on the $L_{x,y,z}$'s.  In this
sense, we will call the $L_{a,b,c}$'s as ``building blocks'' for the regular
representation, or equivalently, ``building blocks'' of $A$.

\section{Comultiplication.  The multiplicative unitary operator}

 We have constructed our $C^*$--algebra $A$ as a strict deformation quantization
of the Poisson--Lie group $\bigl(G,\{\ \}_{\omega}\bigr)$.  We now proceed to
equip $A$ with its quantum group structures.  The first step is to define an 
appropriate comultiplication on it.  An efficient way is to associate a suitable
``multiplicative unitary operator'' \cite{BS}.  That is, we look for a unitary
operator $U$ defined on the Hilbert space ${\mathcal H}\otimes{\mathcal H}$,
such that the ``pentagon equation'' holds (i.\,e. $U_{12}U_{13}U_{23}
=U_{23}U_{12}$) and such that the comultiplication on $A$ is given by
$$
\Delta\phi=U(\phi\otimes1)U^*=U(L_{\phi}\otimes1)U^*,
$$
for $\phi$ contained in the dense subalgebra ${\mathcal A}$ of $A$. 

 To motivate our choice of $U$, let us recall the multiplicative unitary operator
for the ordinary group $C^*$--algebra $C^*(H)$.  It is the operator $V$ defined
on $L^2(H\times H)$ by
\begin{align}
(V\eta)(x,y,z;x',y',z')&=\eta\bigl(x,y,z;(x,y,z)^{-1}(x',y',z')\bigr) 
\notag \\
&=\eta\bigl(x,y,z;x'-x,y'-y,z'-z-\beta(x,y'-y)\bigr). \notag
\end{align}
It is well known \cite{BS}, \cite{ES} that $V$ describes the usual cocommutative
Hopf $C^*$--algebra structure on $C^*(H)$.  Via partial Fourier transform, it may
as well be viewed as an operator on the $(x,y,r)$ variables (still denoted by $V$):
$$
(V\xi)(x,y,r,x',y',r')=\bar{e}\bigl[r'\beta(x,y'-y)\bigr]
\xi(x,y,r+r',x'-x,y'-y,r').
$$

 Since $A$ is essentially a twisted $C^*(H)$, we expect that $V$ needs to be
changed accordingly.  Since the above $V$ represents the regular representation
of $C^*(H)$ \cite{BS}, we expect that the new unitary operator should reflect
the regular representation, $L$, of our twisted group $C^*$--algebra.  So by
using the trick of ``changing of cocycles'' that we mentioned earlier, let us
first consider the following unitary operator $V_{\sigma}$ (where $\sigma$ is
included to emphasize the cocycle) defined on $\mathcal H\otimes\mathcal H$:
$$
(V_{\sigma}\xi)(x,y,r,x',y',r')=\bar{e}\bigl[\eta_{\lambda}(r')\beta
(x,y'-y)\bigr]\xi(x,y,r+r',x'-x,y'-y,r').
$$

 We also have to take into account the point that $A$ should be a quantum version
of $C_{\infty}(G)$.  We will do this by introducing a certain unitary operator
$W$ carrying the information on $G$.  The idea is similar to the ``dual cocycle''
of Landstad \cite{Ld}, \cite{EV}, although $W$ is not exactly a dual cocycle and
$V_{\sigma}$ is not even multiplicative.  Let us consider the following operator
$W$ on $L^2(G\times G)$, motivated by the group multiplication law on $G$:
$$
(W\zeta)(p,q,r,p',q',r')=(e^{\lambda r'})^n\,\zeta(e^{\lambda r'}p,
e^{\lambda r'}q,r,p',q',r').
$$
We may view it as an operator on ${\mathcal H}\otimes{\mathcal H}$, still denoted
by $W$:
$$
(W\xi)(x,y,r,x',y',r')=(e^{-\lambda r'})^n\,\xi(e^{-\lambda r'}x,
e^{-\lambda r'}y,r,x',y',r').
$$
We then incorporate $W$ with $V_{\sigma}$ by defining the unitary operator
$U=WV_{\sigma}$.  We will show in what follows that $U$ is the multiplicative
unitary operator for $A$ we are looking for.  We begin by showing that $U$ is
indeed multiplicative.

\begin{prop}\label{multunitary}
Let $U$ be the unitary operator on $\mathcal H\otimes\mathcal H$ defined by
\begin{align}
U\xi(x,y,r,x',y',r')&=WV_{\sigma}\xi(x,y,r,x',y',r')  \notag \\
&=(e^{-\lambda r'})^n\,\bar{e}\bigl[\eta_{\lambda}(r')\beta(e^
{-\lambda r'}x,y'-e^{-\lambda r'}y)\bigr] \notag \\
&\quad\xi(e^{-\lambda r'}x,e^{-\lambda r'}y,r+r',
x'-e^{-\lambda r'}x,y'-e^{-\lambda r'}y,r').   \notag
\end{align}
Then $U$ is multiplicative.  That is, it satisfies the following ``pentagon
equation'':
$$
U_{12}U_{13}U_{23}=U_{23}U_{12}.
$$
\end{prop}

\begin{proof}
Use Lemma \ref{lem} and calculate:
\begin{align}
(U_{23}&U_{12}\xi)(x_1,y_1,r_1,x_2,y_2,r_2,x_3,y_3,r_3)  \notag \\
=&\,(e^{-\lambda r_3})^n\,\bar{e}\bigl[\eta_{\lambda}(r_3)\beta(e^{-\lambda r_3}
x_2,y_3-e^{-\lambda r_3}y_2)\bigr]  \notag \\
&\,(e^{-\lambda r_2-\lambda r_3})^n\,\bar{e}\bigl[\eta_{\lambda}(r_2+r_3)
\beta(e^{-\lambda(r_2+r_3)}x_1,e^{-\lambda r_3}y_2-e^{-\lambda(r_2+r_3)}y_1)
\bigr]  \notag \\
&\,\xi(e^{-\lambda(r_2+r_3)}x_1,
e^{-\lambda(r_2+r_3)}y_1,r_1+r_2+r_3,e^{-\lambda r_3}x_2-
e^{-\lambda(r_2+r_3)}x_1, \notag \\ 
&\quad e^{-\lambda r_3}y_2-e^{-\lambda(r_2+r_3)}y_1,
r_2+r_3,x_3-e^{-\lambda r_3}x_2,y_3-e^{-\lambda r_3}y_2,r_3)  \notag \\
=&\,(U_{12}U_{13}U_{23}\xi)(x_1,y_1,r_1,x_2,y_2,r_2,x_3,y_3,r_3).   
\notag
\end{align} 
\end{proof}

 For the building block $L_{a,b,c}$ we introduced earlier, define $\Delta L_{a,b,c}$
by
\begin{equation}
\Delta L_{a,b,c}=U(L_{a,b,c}\otimes1)U^*.   \label{(deltadef)}
\end{equation}
Then as an operator on ${\mathcal H}\otimes{\mathcal H}$, we have:
\begin{align}
&(\Delta L_{a,b,c}\xi)(x,y,r,x',y',r') \notag \\
&=\bar{e}\bigl[\eta_{\lambda}(r+r')
\beta(a,e^{-\lambda r'}y-b)+\eta_{\lambda}(r')\beta(a,
y'-e^{-\lambda r'}y)\bigr]\bar{e}\bigl[(r+r')c\bigr] \notag \\
&\quad\xi(x-e^{\lambda r'}a,y-e^{\lambda r'}b,r,
x'-a,y'-b,r').  \label{(deltal)} 
\end{align}
We can show that it is contained in the multiplier algebra $M(A\otimes A)$, which
is rather straightforward (See also the proof of Theorem \ref{comultiplication}
below.).  Moreover, we may regard it as a (continuous) function on $G\times G$
as follows:
\begin{align}
(\Delta L_{a,b,c})(p,q,r,p',q',r')&=\bar{e}\bigl[\langle(e^{\lambda r'}p+p',
e^{\lambda r'}q+q',r+r'),(a,b,c)\rangle\bigr]  \notag \\
&=L_{a,b,c}\bigl((p,q,r)(p',q',r')\bigr) \label{(functiondeltal)}
\end{align}
(Call this function $F\in C_{\infty}(G\times G)$ and use equation \eqref{(rep)}
to compute $(L\otimes L)_F$.  Using (partial) Fourier transform purely formally,
we can show that it agrees with $\Delta L_{a,b,c}$ given by \eqref{(deltal)}).

 In other words, at the level of the building blocks $L_{a,b,c}$, the map $\Delta$
works as the natural comultiplication on $C_{\infty}(G)$.  In view of the fact
that the Poisson structure $\delta$ on $H=G^*$ is linear (see section 1), this
is a desirable choice.  Let us now extend $\Delta$ to the whole algebra $A$ and
obtain our comultiplication:

\begin{theorem}\label{comultiplication}
For $\phi\in{\mathcal A}$, define $\Delta\phi$ by
$$
\Delta\phi=U(\phi\otimes1)U^*=\int_H({\mathcal F}^{-1}\phi)(x,y,z)
\Delta L_{x,y,z}\,dxdydz. 
$$
As before, $\phi$ and  $\Delta\phi$ are actually understood as the operators
$L_{\phi}$ and $(L\otimes L)_{\Delta\phi}$.  Then $\Delta$ can be extended to
a map $\Delta:A\to M(A\otimes A)$, and $\Delta$ is the comultiplication on $A$.
That is, $\Delta$ is a nondegenerate $C^*$--homomorphism satisfying the
coassociativity law: 
$$
(\Delta\otimes\operatorname{id})\Delta\phi
=(\operatorname{id}\otimes\Delta)\Delta\phi.
$$
\end{theorem}

\begin{proof}
It is clear that the formula $\Delta\phi=U(\phi\otimes1)U^*$ defines a
${}^*$--homomorphism, which can be naturally extended to a representation
of $A$ into ${\mathcal B}({\mathcal H}\otimes{\mathcal H})$.

To prove that $\Delta$ carries $A$ into the multiplier algebra $M(A\otimes A)$,
we intend to show that:
$$
(\Delta\phi)(1\otimes g)\in S_{3c}(\Gh/\Gz\times\Gg/\Gq\times\Gh/\Gz\times\Gg/\Gq),
\qquad g\in S_{3c}(\Gh/\Gz\times\Gg/\Gq)\cong{\mathcal A}.
$$
Here, $S_{3c}(\Gh/\Gz\times\Gg/\Gq\times\Gh/\Gz\times\Gg/\Gq)$ is the space of
Schwartz functions having compact support in the $r$ and the $r'$ variables.
This is a dense subspace of $A\otimes A$ (See remark below.).

Let $\xi\in{\mathcal H}\otimes\mathcal H$ and calculate.  We use the change
of variables and the Fourier inversion theorem.  Also, the identity \eqref
{(lemeq)} of Lemma \ref{lem} is very convenient.  We have:
\begin{align}
&(\Delta\phi(1\otimes g)\xi)(x,y,r,x',y',r') \notag \\
&=\int({\mathcal F}^{-1}\phi)(a,b,c)\bar{e}\bigl[\eta_{\lambda}(r+r')
\beta(a,e^{-\lambda r'}y-b)+\eta_{\lambda}(r')\beta(a,y'-e^{-\lambda r'}y)
\bigr] \notag \\
&\qquad\quad\bar{e}\bigl[(r+r')c\bigr]g(\tilde{x},\tilde{y},r')
\bar{e}\bigl[\eta_{\lambda}(r')\beta(\tilde{x},y'-b-\tilde{y})\bigr] \notag \\
&\qquad\quad\xi(x-e^{\lambda r'}a,y-e^{\lambda r'}b,r,x'-a-\tilde{x},
y'-b-\tilde{y},r')\,dadbdcd\tilde{x}d\tilde{y} \notag \\ 
&=(F\xi)(x,y,r,x',y',r')=(L\otimes L)_F\xi(x,y,r,x',y',r'), \notag
\end{align}
where $F$ is defined by
\begin{align} 
F(x,y,r,x',y',r')&=(e^{-2\lambda r'})^n\,\bar{e}\bigl[\eta_{\lambda}(r')
\beta(e^{-\lambda r'}x,y'-e^{-\lambda r'}y)\bigr] \notag \\
&\phi^{\vee}(e^{-\lambda r'}x,e^{-\lambda r'}y,r+r')g(x'
-e^{-\lambda r'}x,y'-e^{-\lambda r'}y,r'). \notag
\end{align}
It is easy to see that $F\in S_{3c}(\Gh/\Gz\times\Gg/\Gq\times\Gh/\Gz\times
\Gg/\Gq)$.  A similar result also holds when we multiply $\Delta\phi$ from
the right.

Since $(\Delta\phi)(1\otimes g)\in A\otimes A$ and $(1\otimes g)(\Delta\phi)
\in A\otimes A$, for an arbitrary $g$ contained in a dense subset of $A$,
we can see that $\Delta\phi\in M(A,A)$, where $M(A,A)=
\{x\in M(A\otimes A):x(1_{M(A)}\otimes A)+(1_{M(A)}\otimes A)x\in A\otimes A\}$.
It is customary to require (see \cite{Va}, \cite{BS}) that the comultiplication
$\Delta$ takes values in $M(A,A)$.  This is done so that one is able to discuss
the notion of ``left invariant'' Haar weight on $A$.  

Actually, we can improve the statement even further by observing that $(\Delta
\phi)(1\otimes g)$'s form a total set (with respect to the $L^1$--norm) in the
Schwartz space $S_{3c}(\Gh/\Gz\times\Gg/\Gq\times\Gh/\Gz\times\Gg/\Gq)$.  We may
check this using the expression given above.  Since the Schwartz space is in
turn a dense subspace of $A\otimes A$, this is enough to show that $\Delta$ is
also non-degenerate. 

 Finally, the coassociativity of $\Delta$ follows from the fact that $U$ is
multiplicative.  For $\phi\in{\mathcal A}$, we would have:
$$
U_{12}U_{13}(\phi\otimes1\otimes1)U_{13}^*U_{12}^*
=U_{23}U_{12}(\phi\otimes1\otimes1)U_{12}^*U_{23}^*.
$$
But by definition of $\Delta$, this is just:
$$
(\Delta\otimes\operatorname{id})\Delta\phi
=(\operatorname{id}\otimes\Delta)\Delta\phi.
$$
\end{proof}

\begin{rem}
For the above proof to be complete, we need to show that the Schwartz space
$S_{3c}(\Gh/\Gz\times\Gg/\Gq\times\Gh/\Gz\times\Gg/\Gq)$ is a dense subset
of $A\otimes A$.  So consider the natural injection from $S_{3c}(\Gh/\Gz\times
\Gg/\Gq\times\Gh/\Gz\times\Gg/\Gq)$ into ${\mathcal B}({\mathcal H}\otimes
{\mathcal H})$, which is continuous with respect to the $L^1$--norm on the
Schwartz space and the $C^*$--norm on ${\mathcal B}({\mathcal H}\otimes
{\mathcal H})$.  Under this natural injection, the algebraic tensor product
$S_{3c}(\Gh/\Gz\times\Gg/\Gq)\odot S_{3c}(\Gh/\Gz\times\Gg/\Gq)$ is sent into
a dense subset of the algebraic tensor product $A\odot A$.  Since any element
in $S_{3c}(\Gh/\Gz\times\Gg/\Gq\times\Gh/\Gz\times\Gg/\Gq)$ can be approximated
by elements of $S_{3c}(\Gh/\Gz\times\Gg/\Gq)\odot S_{3c}(\Gh/\Gz\times\Gg/\Gq)$
in the $L^1$--norm, we conclude that $S_{3c}(\Gh/\Gz\times\Gg/\Gq\times\Gh/\Gz
\times\Gg/\Gq)$ is mapped into a dense subset of $A\otimes A$.
\end{rem}

  Our choice of $\Delta L_{a,b,c}$ (equation \eqref{(deltal)}) together with
the above theorem means that the comultiplication remains the ``same'' while
the algebra is being deformed (i.\,e. Our deformation is a {\em preferred\/}
deformation \cite{G1}, \cite{G2}.).  In this way, we obtained our Hopf
$C^*$--algebra $(A,\Delta)$.

\begin{defn}
(\cite{Va}, \cite{BS})
By a {\em Hopf $C^*$--algebra\/}, we mean a pair $(B\Delta)$, where $B$ is
a $C^*$--algebra and $\Delta$ is a comultiplication (satisfying the conditions
given in Theorem \ref{comultiplication}).
\end{defn}

 It may not be evident, but our construction is closely related with Baaj and
Skandalis's construction of Hopf $C^*$--algebras via ``matched pair'' ({\em
couple assorti\/}) and ``bicrossed product'' ({\em biproduit crois\'e\/}) of
``Kac Systems'' \cite[\S8]{BS} (Similar work at the algebraic level is done by
Majid \cite{Mj}.).

 To be a little more specific, the abelian groups $H/Z$ and $(\Gg/\Gq,+)$ (or
in terms of Hopf $C^*$--algebras, $C^*(H/Z)$ and $C_{\infty}(\Gg/\Gq)$) form
a matched pair.  From this, we can form a ``twisted'' bicrossed product, using
the notion of cocycles satisfying certain equivariance condition.  The
multiplicative unitary operator associated with this ``matched pair with
cocycle'' construction is {\em regular\/} \cite{BS}.  Although our construction
of $(A,\Delta)$ and Baaj and Skandalis' method are rather different, we can
still show that our multiplicative unitary operator $U$ for $A$ coincides with
the multiplicative unitary operator for this twisted bicrossed product.

 We do not intend to prove the regularity of $U$ directly (However, the result
in the proof of Theorem \ref{comultiplication} that $(\Delta\phi)(1\otimes g)$'s
form a total set in $A\otimes A$ is very much related.).  Instead, let us refer
to the above discussions and summarize the result in the following:

\begin{prop}
Let $U$ be the unitary operator defined as in Proposition \ref{multunitary}.
It is a ``regular'' multiplicative unitary operator, in the sense of Baaj and
Skandalis. 
\end{prop}

 This result gives our construction an axiomatically sound basis:  If we start
from the multiplicative unitary operator $U$, its associated Hopf $C^*$--algebra
is exactly $(A,\Delta)$.  Also associated with the regular multiplicative unitary
operator is the dual Hopf $C^*$--algebra $(\hat{A},\hat{\Delta})$.  In our case, 
$\hat{A}$ is essentially the group $C^*$--algebra $C^*(G)$.  And $\hat{A}\cong
C^*(G)$ is a deformation quantization of $H$, equipped with the (linear) Poisson
bracket given by $\delta$ defined in section 1.  In this way, we see that the
duality between $H$ and $G$ as Poisson--Lie groups corresponds nicely to the Hopf
$C^*$--algebra duality between $\hat{A}$ and $A$, in terms of the multiplicative
unitary operator $U$.

\bigskip

\begin{center} {\sc Appendix: Deformation quantization of $\tilde{G}$}
\end{center}

\smallskip

 Recall from section 1 that the Lie bialgebra structure $\delta$ on $\Gh$ actually
came from the Lie bialgebra structure $\tilde{\delta}$ on the extended Heisenberg
Lie algebra $\tilde{\Gh}$.  So far, we only considered the deformation quantization
of $\bigl(G,\{\ \}_{\omega}\bigr)$, which is the dual Poisson--Lie group of the
(nilpotent) Poisson--Lie group $H$ or Lie bialgebra $(\Gh,\delta)$.  We have been
avoiding the discussion of $\tilde{H}$ and its dual Lie group $\tilde{G}$, because
$\tilde{H}$ is not nilpotent.

 Usually, there are some technical difficulties to correctly formulate the notion
of ``strict'' deformation quantization of $C_{\infty}(\tilde{G})$, if $\tilde{H}$
is not nilpotent.  Some modifications of the ``strictness condition'' should be
necessary.  See \cite{BJKp1}, \cite{Rf3}.  But in our case, if we are willing to
compromise a little on shrinking the space on which the deformed multiplication is
defined, we are still able to find a quantum version of $C_{\infty}(\tilde{G})$,
with the aid of multiplicative unitary operators.  We are going to define below a
multiplicative unitary operator $\tilde{U}$, using the trick of ``changing of
cocycles'' as before.  The multiplicative unitary we obtain will again be regular.

 By \cite{BS}, given a regular multiplicative unitary $\tilde{U}\in{\mathcal B}
(\tilde{\mathcal H}\otimes\tilde{\mathcal H})$, there corresponds an algebra
${\mathcal A}(\tilde{U})\subseteq{\mathcal B}(\tilde{\mathcal H})$ such that
its norm closure gives a $C^*$--algebra $\tilde{A}$.  Usually, ${\mathcal A}
(\tilde{U})$ is kind of an $L^1$--algebra.  In our case, it will be the twisted
group algebra whose twisted convolution is given by the cocycle associated to
the definition of $\tilde{U}$.  Since we prefer to have our multiplication defined
at the level of continuous functions on $\tilde{G}$, we will consider a certain
subspace $\tilde{\mathcal A}$ of $S(\tilde{G})$, to express our multiplication.

 The following construction is indeed a deformation quantization of $\tilde{G}$.
The verification of this will be left to the reader.

\begin{ex}\label{ex}
Let $\tilde{H}$ be the extended Heisenberg Lie group with the group law
defined by
\small
$$
(x,y,z,w)(x',y',z',w')=\bigl(x+e^{w}x',y+e^{-w}y',z+z'+(e^{-w})\beta(x,y'),
w+w'\bigr).
$$
\normalsize
This is clearly the Lie group corresponding to the extended Heisenberg Lie algebra
$\tilde{\Gh}$ defined in section 1.  We use the $w$ variable here to express the
vectors in $\operatorname{span}(\mathbf{d})$.  Consider the dual Poisson--Lie group
$\tilde{G}$ of $\tilde{H}$ defined by the multiplication:
\small
$$
(p,q,r,s)(p',q',r',s')=(e^{\lambda r'}p+p',e^{\lambda r'}q+q',r+r',
s+s').
$$
\normalsize
It is easy to see that the above $\tilde{G}$ is indeed the Lie group associated
with the Lie algebra $\tilde{\Gg}$ defined in Corollary of Proposition
\ref{bialgebra}. To describe its deformation quantization, it is convenient
to work in the space of $(x,y,r,w)$ variables, $S(\tilde{\Gh}/\Gz\times
\tilde{\Gg}/\tilde{\Gq})$.  Here $\tilde{\Gg}=\tilde{\Gh}^*$ and $\tilde{\Gq}
=\Gz^{\bot}$ in $\tilde{\Gg}$.

\smallskip
\noindent (multiplication): Consider the subspace of $S(\tilde{\Gh}/\Gz\times
\tilde{\Gg}/\tilde{\Gq})$ having compact support in both the $r$ and the $w$
variables.  Let $\tilde{\mathcal A}$ be its image in $S(\tilde{G})$ under partial
Fourier transform in the $(x,y,w)$ variables, still denoted by ${}^{\wedge}$.
On $\tilde{\mathcal A}$ we define the deformed multiplication by
\small
\begin{align}
(\phi&\times\psi)(p,q,r,s)  \notag\\
&=\int\phi^{\vee}(\tilde{x},\tilde{y},r,\tilde{w})\psi^{\vee}
(e^{-\tilde{w}}x-e^{-\tilde{w}}\tilde{x},e^{\tilde{w}}y
-e^{\tilde{w}}\tilde{y},r,w-\tilde{w})  \notag
\\ &\qquad\quad\bar{e}\bigl[\eta_{\lambda}(r)\beta(\tilde{x},
y-\tilde{y})\bigr]\bar{e}[p\cdot x+q\cdot y+sw]\,d\tilde{x}d\tilde{y}
d\tilde{w}dxdydw,  \notag
\end{align}
\normalsize
where ${}^{\vee}$ is the (inverse) partial Fourier transform in the $(p,q,s)$
variables.  This definition of $\times$ is motivated by the fact that the
Poisson bracket on $\tilde{G}$ is essentially the extension of the linear
Poisson bracket on $(\tilde{\Gh}/\Gz)^*$ by a cocycle.  We follow the method
of \cite{BJKp1}.  Our $C^*$--algebra $\tilde{A}$ will then be defined as the
enveloping $C^*$--algebra of $\bigl(\tilde{\mathcal A},\times\bigr)$.

\smallskip
\noindent (comultiplication): Define the following unitary operators on
$\tilde{\mathcal H}\otimes\tilde{\mathcal H}$, where $\tilde{\mathcal H}$
is the space of $L^2$--functions on the $(x,y,r,w)$ variables.
\small
\begin{align}
\tilde{W}\xi(x,y,r,w&,x',y',r',w')=(e^{-\lambda r'})^n\,
\xi(e^{-\lambda r'}x,e^{-\lambda r'}y,r,w,x',y',r',w')  \notag \\
\tilde{V}_{\sigma}\xi(x,y,r,w&,x',y',r',w')=\bar{e}\bigl[\eta_{\lambda}
(r')\beta(x,y'-y)\bigr]  \notag \\
&\quad\xi(x,y,r+r',w,e^{-w}x'-e^{-w}x,e^{w}y'-e^{w}y,r',w'-w).
\notag
\end{align}
\normalsize
Let $\tilde{U}=\tilde{W}\tilde{V}_{\sigma}$.  Then we would have:
\small
\begin{align}
\tilde{U}&\xi(x,y,r,w,x',y',r',w')=(e^{-\lambda r'})^n\,\bar{e}\bigl[
\eta_{\lambda}(r')\beta(e^{-\lambda r'}x,y'-e^{-\lambda r'}y)\bigr]  
\notag \\
&\xi(e^{-\lambda r'}x,e^{-\lambda r'}y,r+r',w,e^{-w}x'
-e^{-\lambda r'-w}x,e^{w}y'-e^{-\lambda r'+w}y,r',w'-w).  \notag 
\end{align}
\normalsize
It is again a multiplicative unitary operator.  Thus we may define the
comultiplication on $\tilde{A}$ by $\tilde{\Delta}\phi=\tilde{U}(\phi\otimes1)
\tilde{U}^*$.  Since it will be useful in later calculations, let us write
down the explicit formula for the comultiplication of the building block
$L_{a,b,c,d}$, for $(a,b,c,d)\in\tilde{H}$.

 For $(a,b,c,d)\in\tilde{H}$, the building block $L_{a,b,c,d}$ is the
operator on $\tilde{\mathcal H}$ defined similarly as in equation
\eqref{(block)} earlier:
\small
$$
(L_{a,b,c,d}\xi)(x,y,r,w)=\bar{e}(rc)\bar{e}\bigl[\eta_{\lambda}(r)
\beta(a,y-b)\bigr]\xi(e^{-d}x-e^{-d}a,e^d y-e^d b,r,w-d).
$$
\normalsize
So $\tilde{\Delta}L_{a,b,c,d}=\tilde{U}(L_{a,b,c,d}\otimes1)\tilde{U}^*$
is an operator on $\tilde{\mathcal H}\otimes\tilde{\mathcal H}$ defined by
\small
\begin{align}
&(\tilde{\Delta}L_{a,b,c,d}\xi)(x,y,r,w,x',y',r',w')  \notag \\
&=\bar{e}\bigl[\eta_{\lambda}(r)\beta(e^{\lambda r'}a,y-e^{\lambda r'}b)+
\eta_{\lambda}(r')\beta(a,y'-b)\bigr]\bar{e}\bigl[(r+r')c\bigr] \notag \\
&\quad\xi(e^{-d}x-e^{\lambda r'-d}a,e^{d}y-e^{\lambda r'+d}b,r,w-d,
e^{-d}x'-e^{-d}a,e^{d}y'-e^{d}b,r',w'-d).  \notag
\end{align}
\normalsize
\end{ex}

\section{Counit and antipode}

 We return to the construction of the remaining quantum group structures
for our Hopf $C^*$--algebra $(A,\Delta)$.  Similar results will hold for
$(\tilde{A},\tilde{\Delta})$ since we only need to change the groups
accordingly and use the appropriate cocycles.  So in this section and the
next, we will exclusively study about our specific example $(A,\Delta)$.
Since $A$ is our candidate for the ``quantum $C_{\infty}(G)$'', we expect
that its quantum group structures will come from the corresponding group
structures on $G$.

 First, the choice for the counit is rather obvious:
\begin{theorem}\label{counit}
There exists a unique continuous linear map $\epsilon:A\to\mathbb{C}$ such that
$$
\epsilon(\phi)=\phi(0,0,0),
$$
for $\phi\in{\mathcal A}=S_{3c}(G)$.  Then $\epsilon$ is a counit for 
$(A,\Delta)$.  That is, $\epsilon$ is a $C^*$--homomorphism from $A$
into $\mathbb{C}$ satisfying the condition:
$$
(\operatorname{id}\otimes\epsilon)\Delta=(\epsilon\otimes
\operatorname{id})\Delta=\operatorname{id}.
$$
\end{theorem}

\begin{proof}
For $\phi\in{\mathcal A}$,
\begin{align}
\epsilon(\phi)=\phi(0,0,0)&=\int({\mathcal F}^{-1}\phi)(x,y,z)\,dxdydz 
\notag \\
&=\int({\mathcal F}^{-1}\phi)(x,y,z)\epsilon(L_{x,y,z})\,dxdydz,   
\notag
\end{align}
where $\epsilon(L_{x,y,z})\equiv1$.  In other words, $\epsilon$ is actually
the trivial representation of ${\mathcal A}$.  On the other hand, we may write
$\epsilon(\phi)$ as:
$$
\epsilon(\phi)=\phi(0,0,0)=\int\phi^{\vee}(x,y,0)\,dxdy,
$$
which shows that $\epsilon$ is continuous with respect to the $L^1$-norm.
So $\epsilon$ has a continuous linear extension to the $L^1$--algebra 
$L^1\bigl(H/Z,C_{\infty}(\Gg/\Gq),\sigma\bigr)$.  But since we have already
seen that $\epsilon$ is a ${}^*$--representation on ${\mathcal A}$, this
extension is also a ${}^*$--representation.  Therefore, it can be further
extended to a ${}^*$--representation on $A\cong C^*\bigl(H/Z,C_{\infty}(\Gg/\Gq),
\sigma\bigr)$.

 Next, let us prove the equality for our building block $\Delta L_{x,y,z}
\in M(A\otimes A)$.  Using the realization of $\Delta L_{x,y,z}$ as a continuous
function on $G\times G$ (equation \eqref{(functiondeltal)}), we have:
$$
(\operatorname{id}\otimes\epsilon)\Delta L_{x,y,z}(p,q,r)=\bar{e}
\bigl[\langle(p,q,r),(x,y,z)\rangle\bigr]=L_{x,y,z}(p,q,r),
$$
and similarly for the other half of the equality.  By the 
definition of $\Delta$, we have proved that: 
$$
(\operatorname{id}\otimes\epsilon)(\Delta\phi)=\phi=(\epsilon\otimes
\operatorname{id})(\Delta\phi),\qquad\phi\in{\mathcal A}.
$$
\end{proof}

 The antipode (or coinverse) is usually defined as an anti-automorphism \cite{Va},
\cite{Wr7}.  Let us follow the method which has been used by several authors,
beginning as early as the work by Kac and Paljutkin \cite{KP}.

 Consider the operation $\dagger$ on $\mathcal A$ defined by
$$
\phi^{\dagger}(p,q,r)=\overline{\phi(-e^{-\lambda r}p,
-e^{-\lambda r}q,-r)}.
$$
The bar means the complex conjugation.  Then define $\kappa:{\mathcal A}
\to{\mathcal A}$ by
$$
\kappa(\phi)=(\phi^*)^{\dagger}=(\phi^{\dagger})^* ,
$$
where $\phi^*$ is the $C^*$--involution defined in Proposition \ref{invol}.
Explicitly, we have:
\begin{align}
\kappa(\phi)(p,q,r)&=\int\phi(-e^{-\lambda r}\tilde{p},-e^{-\lambda r}
\tilde{q},-r)\bar{e}\bigl[(p-\tilde{p})\cdot x+(q-\tilde{q})
\cdot y\bigr]  \notag \\
&\quad\quad\quad\bar{e}\bigl[\eta_{\lambda}(r)\beta(x,y)\bigr]
\,d\tilde{p}d\tilde{q}dxdy.  \label{(kappa)}
\end{align} 
In the commutative case (i.\,e. $\beta\equiv0$), this is none other than:
$$
\kappa(\phi)(p,q,r)=\phi(-e^{-\lambda r}p,-e^{-\lambda r}q,-r)
=\phi\bigl((p,q,r)^{-1}\bigr),
$$
which is just taking the inverse in $G$.

 Let us now try to define $\kappa$ at the operator level.  Motivated by the
operation $\dagger$ above, we first define an involutive operator $T$ on
$\mathcal H$ by
$$
T\xi(x,y,r)=(e^{\lambda r})^n\,\overline{\xi(e^{\lambda r}x,
e^{\lambda r}y,-r)}.
$$

\begin{lem}
Let $T$ be the operator defined above.  Then $T$ is conjugate 
linear, isometric, and involutive (i.\,e. $T^2=1$.).  Moreover,
$$
T\phi T=\phi^{\dagger},
$$
where $\phi,\phi^{\dagger}\in{\mathcal A}$ are viewed as operators.  
We thus have $TAT=A$.
\end{lem}

\begin{proof}
We will just verify the equation $T\phi T=\phi^{\dagger}$.  The other
assertions are straightforward.  We have:
\begin{align}
(T\phi T\xi)(x,y,r)&=(e^{\lambda r})^n\,\overline{(\phi T\xi)
(e^{\lambda r}x,e^{\lambda r}y,-r)} \notag \\
&=\int\overline{\phi(p,q,-r)}\bar{e}(p\cdot x'+q\cdot y')
e\bigl[\eta_{\lambda}(-r)\beta(x',e^{\lambda r}y-y')\bigr]  \notag \\
&\qquad\quad\xi(x-e^{-\lambda r}x',y-e^{-\lambda r}y',r)\,dpdqdx'dy'  
\notag \\
&=\phi^{\dagger}\xi(x,y,r). \notag
\end{align}
\end{proof} 

\begin{prop}\label{kappa}
Let the map $\kappa:A\to A$ be defined by $\kappa(\phi)=T\phi^*T$, for 
$\phi\in{\mathcal A}$.  Then $\kappa$ is an anti-automorphism on $A$.
At the function level, $\kappa(\phi)$ agrees with equation \eqref{(kappa)}.
Moreover, $\kappa$ satisfies the condition:
$$
(\kappa\otimes\kappa)(\Delta\phi)=\Sigma\bigl(\Delta(\kappa\phi)\bigr),
$$
where $\Sigma:A\otimes A\to A\otimes A$ denotes the flip.
\end{prop}

\begin{proof}
The proof that $\kappa$ is an anti-automorphism follows immediately 
from the previous lemma.  Since $\kappa(\phi)=(\phi^{\dagger})^*$ on the
functions, to prove the last condition we only need to check the 
following equation:
$$
(T\otimes T)(\Delta\phi)(T\otimes T)\xi=(\Sigma(\Delta\phi^
{\dagger}))(\Sigma\xi),\qquad\xi\in{\mathcal H}\otimes{\mathcal H}.
$$
Here $\Sigma$ also denotes the flip on ${\mathcal H}\otimes{\mathcal H}$.  The
calculation is straightforward.
\end{proof}

 In this way, we showed that $(A,\Delta,\epsilon,\kappa)$ is a {\em counital,
coinvolutive Hopf $C^*$--algebra\/} in the sense of \cite{Va}.  However, some
more axioms are needed to make the map $\kappa$ to be reasonably considered as the
antipode.  For instance, in the purely algebraic setting of Hopf algebra theory
\cite{S}, \cite{M}, the requirement for the antipode is  given by the following
equation:
\begin{equation}\label{(antipode)}
m(\operatorname{id}\otimes\kappa)\Delta=m(\kappa\otimes\operatorname
{id})\Delta = \epsilon(\cdot)1,
\end{equation}
where $m:A\otimes A\to A$ is the multiplication.  

 In the operator algebra setting, the multiplication map $m$ is not continuous
for the operator norms in general.  Because of this, we approach a little
differently rather than just translating the above formulation.  Motivated by
Kac algebra theory, the antipode is usually dicussed together with the notion
of the Haar weight.  See Proposition \ref{haar} in the next section.

 Nevertheless, at least at the function level, the algebraic condition \eqref
{(antipode)} can be readily verified for our $(A,\Delta,\epsilon,\kappa)$.
The calculation of this claim is as follows.  This will give us some modest
justification for our particular definition of $\kappa$.

 Using the definition of $\Delta$ and the fact that $\Delta L_{a,b,c}$ can be
regarded as a continuous function on $G$, we have for $\phi\in{\mathcal A}$,
\begin{align}
(\operatorname{id}&\otimes\kappa)\Delta\phi(p,q,r,p',q',r')  \notag \\
&=\int\bar{e}\bigl[(p'-\tilde{p})\cdot x+(q'-\tilde{q})\cdot y\bigr]
\bar{e}\bigl[\eta_{\lambda}(r')\beta(x,y)\bigr]
({\mathcal F}^{-1}\phi)(a,b,c)  \notag \\
&\quad\quad\quad\bar{e}\bigl[(e^{-\lambda r'}p-
e^{-\lambda r'}\tilde{p})\cdot a+(e^{-\lambda r'}q-e^{-\lambda r'}
\tilde{q})\cdot b+(r-r')c\bigr]   \notag\\
&\quad\quad\quad d\tilde{p}d\tilde{q}dxdydadbdc \notag \\
&=\int(\psi^1_{a,b,c}\otimes\psi^2_{a,b,c})(p,q,r,p',q',r')\,dadbdc.
\notag
\end{align}
Here for a fixed $(a,b,c)\in H$, 
\begin{align}
\psi^1_{a,b,c}(p,q,r)&=\bar{e}(p\cdot a+q\cdot b+rc) \notag \\
\psi^2_{a,b,c}(p',q',r')&=(e^{2\lambda r'})\,e(p'\cdot a+q'\cdot b+r'c)
\bar{e}\bigl[\eta_{\lambda}(r')\beta(a,b)\bigr].  \notag  \\
&\quad ({\mathcal F}^{-1}\phi)(e^{\lambda r'}a,e^{\lambda r'}b,c).
\notag
\end{align}
Since we have:
$$
(\psi^1_{a,b,c}\times\psi^2_{a,b,c})(p,q,r)=(e^{2\lambda r})
({\mathcal F}^{-1}\phi)(e^{\lambda r}a,e^{\lambda r}b,c),
$$
it follows that:
$$
m\bigl((\operatorname{id}\otimes\kappa)\Delta\phi\bigr)(p,q,r)=\int
({\mathcal F}^{-1}\phi)(a,b,c)\,dadbdc=\phi(0,0,0)=\epsilon(\phi)1.
$$
Similarly, we can also verify: $m\bigl((\kappa\otimes\operatorname{id})
\Delta\phi\bigr)=\epsilon(\phi)1$.

\section{Haar weight}

 Since the group law on $G$ has been chosen such that Lebesgue measure $dpdqdr$
on the underlying vector space is its Haar measure, we expect more or less
the same in the quantum case.  So let us define the linear functional $h$
on $\mathcal A$ by
\begin{equation}\label{(haarfunctional)}
h(\phi)=\int\phi(p,q,r)\,dpdqdr.
\end{equation}
We intend to show that $h$ is the appropriate Haar weight on our Hopf
$C^*$--algebra $A$.

 Ideally, the definition of locally compact quantum groups would be formulated
so that the existence of Haar weights follows only from the definition.  At
present, the definition of Haar weight and its left invariance property are
not completely agreed upon and the existence of Haar weight has to be assumed
in the definition of quantum groups.  In particular, the definition of the
antipode is closely tied to that of the Haar weight.  See \cite{MN}, \cite{MNW},
\cite{Wr7}, \cite{KuVD}, \cite{KuV}.

 Because of this, instead of trying to be very rigorous, we plan to give only a
reasonable justification of our choice for $h$.  What we do in the following is
immitating the theory of Kac algebras \cite{ES}.

 Since $h$ is well-defined at the level of a dense subspace of functions (i.\,e.
in ${\mathcal A}$), it is a densely defined weight on $A$.  As we see in the next
proposition, it is actually a faithful trace.

\begin{prop}
Let $h$ be defined on ${\mathcal A}$ by equation \eqref{(haarfunctional)}.  Then
$h$ is a faithful trace.
\end{prop}

\begin{proof}
Let $\phi\in{\mathcal A}$.  Then by using change of variables and Fourier inversion
theorem, we have:
\begin{align}
h(\phi^*\times\phi)&=\int\bar{e}\bigl[(p-p')\cdot\tilde{x}\bigr]\phi^*(p,q,r)
\phi\bigl(p,q+\eta_{\lambda}(r)\tilde{x},r\bigr)\,dp'd\tilde{x}  \notag\\
&=\int\overline{\phi(p,q,r)}\phi(p,q,r)\,dpdqdr=\|\phi\|_2^2,  \notag
\end{align}
and similarly,
$$
h(\phi\times\phi^*)=\int\phi(p,q,r)\overline{\phi(p,q,r)}\,dpdqdr=\|\phi\|_2^2.
$$
Here, $\phi^*$ is the $C^*$--involution given in Proposition \ref{invol}.  From
these equations, we can see that $h$ is a faithful trace.
\end{proof}

 To correctly define the Haar weight, we have to further require some ``lower
semi-continuity condition'' (corresponding to the notion of normal weights in
von Neumann algebra setting, like Kac algebras) and ``semi-finiteness'', as well
as the ``left invariance property''.  Since this will make our discussion very
technical, let us overlook the details and give only a brief discussion on the
left invariance property of $h$.

\begin{prop}\label{haar}
For $\phi,\psi\in{\mathcal A}$, the weight $h$ satisfies the following left
invariance property:
\begin{equation}\label{(haar)}
(\operatorname{id}\otimes h)\bigl((1\otimes\phi)(\Delta\psi)\bigr)
=\kappa\bigl((\operatorname{id}\otimes h)((\Delta\phi)(1\otimes\psi))
\bigr),
\end{equation}
where $\kappa$ is the (antipodal) map defined in Proposition \ref{kappa}.
\end{prop}

\begin{proof}
Even for $\phi,\psi\in{\mathcal A}$, the expressions $(1\otimes\phi)(\Delta\psi)$
and $(\Delta\phi)(1\otimes\psi)$ do not necessarily belong to the algebraic
tensor product ${\mathcal A}\odot{\mathcal A}$ (See the proof of Theorem
\ref{comultiplication}, where we calculated $(\Delta\phi)(1\otimes\psi)$.).
Therefore, for the left and right sides of the equation \eqref{(haar)} to make
sense, $\operatorname{id}\otimes h$ has to be defined more carefully.

This extension can be done using the notion of operator valued weights
\cite{Haag}.  But unlike in \cite{Haag}, since we are dealing with
$C^*$--algebra weights (\cite{Cm}), we have to modify the definitions
accordingly.  In short, we regard $\operatorname{id}\otimes h$ as the
tensor product of two faithful, semi-finite, lower semi-continuous operator
valued weights, on $A\otimes A$ having values in $A\otimes\mathbb{C}\cong A$.
To be able to define this more rigorously, there are efforts being made
introducing somewhat stronger condition of lower semi-continuity \cite{QV1},
\cite{Ku}.

In our case, since we know that $(1\otimes\phi)(\Delta\psi)$ and $(\Delta\phi)
(1\otimes\psi)$ are contained in $S_{3c}(\Gh/\Gz\times\Gg/\Gq\times\Gh/\Gz
\times\Gg/\Gq)$ and since the elements in this Schwartz space can be approximated
by elements in $S_{3c}(\Gh/\Gz\times\Gg/\Gq)\odot S_{3c}(\Gh/\Gz\times\Gg/\Gq)$,
we know how to define $(\operatorname{id}\otimes h)\bigl((1\otimes\phi)
(\Delta\psi)\bigr)$ and $(\operatorname{id}\otimes h)\bigl((\Delta\phi)
(1\otimes\psi)\bigr)$ under the extension.  So let us set aside the aforementioned
technical details and try to verify the above equation.  Through long but
elementary calculations, we obtain:
\begin{align}
&(\operatorname{id}\otimes h)\bigl((1\otimes\phi)(\Delta\psi)\bigr)
(p,q,r)  \notag\\
&\quad=\int \bar{e}\bigl[(e^{\lambda r'}p+\tilde{p}-\tilde{\tilde{p}})
\cdot a+(e^{\lambda r'}q+\tilde{q}-\tilde{\tilde{q}})\cdot b\bigr]
e\bigl[\eta_{\lambda}(r')\beta(a,b)\bigr] \notag \\
&\quad\quad\quad\quad\phi(\tilde{p},\tilde{q},r')\psi(\tilde{\tilde{p}},
\tilde{\tilde{q}},r+r')\,d\tilde{p}d\tilde{q}d\tilde{\tilde{p}}
d\tilde{\tilde{q}}dadbdr' \notag \\
&\quad=\kappa\bigl((\operatorname{id}\otimes h)((\Delta\phi)(1\otimes\psi))
\bigr)(p,q,r)  \notag
\end{align}
for $\phi$ and $\psi$ in ${\mathcal A}$.
\end{proof}

 In the commutative case, equation \eqref{(haar)} is none other than
$$
\int_G\phi(g')\psi(gg')\,d\mu(g')=\int_G\phi(g^{-1}g')\psi(g')\,d\mu(g'),
$$
which exactly describes the left invariance condition.  Actually, equation 
\eqref{(haar)} is the defining condition for the Haar weight in Kac algebra
theory \cite{Va}, \cite{ES}.

 It is true that there are still some technical details to be taken care of.
Having said this, we may conclude from Proposition \ref{haar} that $h$ is the
appropriate {\em haar weight\/} for $A$.  Also from the proposition, we can
say that the map $\kappa$ we have been using is a legitimate {\em antipode\/}
for $A$.

 Thus our Hopf $C^*$--algebra $(A,\Delta,\epsilon,\kappa)$ together with the
Haar weight $h$ on it can be regarded as a {\em locally compact quantum group\/}.
Although we did not give the precise definition of general locally compact
quantum groups, any reasonable definition should allow our specific example
as a special case.

 Meanwhile, since our group $G$ is not unimodular, we expect that our Haar weight
should also carry certain non-unimodularity properties.  One such is given below:

\begin{prop}
The Haar weight $h$ is not invariant under the antipode $\kappa$.  That is,
there exists $\phi\in{\mathcal A}$ such that
$$
h\bigl(\kappa(\phi)\bigr)\ne h(\phi).
$$
\end{prop}

\begin{proof}
Since
$$
h\bigl(\kappa(\phi)\bigr)=\int\phi(-e^{-\lambda r}p,-e^{-\lambda r}q,-r)
\,dpdqdr,
$$
it is clear that we have $h\bigl(\kappa(\phi)\bigr)\ne h(\phi)$, in general.
\end{proof}

 It is noteworthy that we have a non-unimodular Haar weight as opposed to many
other examples \cite{Rf5}, \cite{SZ}, \cite{VD}, \cite{Ld}.  It will be 
interesting to study its consequences and properties more thoroughly, 
especially in relation to the duality theory.  For the time being, however,
we will leave this as a future project.

 As a final remark, we point out that the regular representation $L$ we have
been using is essentially the GNS representation with respect to $h$ (which is
a faithful trace).  The partial Fourier transform provides the equivalence.
This observation displays the importance the Haar weight has in both theory
and construction of locally compact quantum groups.

\section{Quantum universal $R$--matrix}

 For the QUE algebra counterparts for our Hopf $C^*$--algebra $(A,\Delta)$
(for instance, $U_{\hbar}(\Gh)$ in \cite{ALT} or $H(1)_q$ in \cite{CGST}),
the so-called ``quantum universal $R$--matrix'' have been successfully
constructed.  In our case also, once we modify the definition of the
universal $R$--matrix so that it is consistent with our $C^*$--algebra
language, we can do the same.

 Our definition given below is essentially the same one used in the QUE
algebra or more general Hopf algebra setting (See \cite{Dr}, \cite{CP}.).
Note that we require the $R$--matrix to be contained in a multiplier algebra
(This is consistent with the definition of the comultiplication, which is a
multiplier algebra valued map.).  Since any nondegenerate representation
of a $C^*$--algebra can be uniquely extended to its multiplier algebra,
an element being in a multiplier algebra also means that it has an image
under any representation of the $C^*$--algebra.  

\begin{defn}\label{alcomm}
Let $(B,\Delta)$ be a Hopf $C^*$--algebra, where $\Delta$ is its
comultiplication.  We will say that $B$ is {\em almost cocommutative\/},
if there exists an invertible element $R\in M(B\otimes B)$ such that
\begin{equation}
(\Sigma\circ\Delta)(\phi)=R\Delta(\phi)R^{-1},\quad \phi\in B
\end{equation}
where $\Sigma$ is the flip.  We will denote the {\em opposite
comultiplication\/} by $\Delta^{\operatorname{op}}=\Sigma\circ\Delta$.
\end{defn}

 The element $R$ above cannot be arbitrary, since the opposite comultiplication
$\Delta^{\operatorname{op}}$ should also be coassociative.  The following
condition, though a little stronger than is needed to assure the coassociativity
of $\Delta^{\operatorname{op}}$, defines the quantum universal $R$--matrix. 

\begin{defn}\label{qtrR}
An almost cocommutative Hopf $C^*$--algebra $(B,R)$ is said to be {\em
quasitriangular\/}, if $R$ satisfies the so-called quantum Yang--Baxter
equation (QYBE): $R_{12}R_{13}R_{23}=R_{23}R_{13}R_{12}$,
and also satisfies:  
\begin{equation}\label{(qt)}
(\Delta\otimes\operatorname{id})(R)=R_{13}R_{23}\quad {\text {and}}
\quad(\operatorname{id}\otimes\Delta)(R)=R_{13}R_{12}.
\end{equation}
It is called {\em triangular\/}, if it is quasitriangular and in addition,
$R_{21}=R^{-1}$.  If $B$ is quasitriangular, such an element $R$ will be
called a {\em quantum universal R--matrix\/}.
\end{defn}

 If $R$ satisfies equation \eqref{(qt)}, the QYBE for $R$ automatically
follows from the coassociativity of $\Delta^{\operatorname{op}}$ \cite{CP}. 
The QYBE is a quantum version of the classical Yang--Baxter equation (CYBE)
\cite{Dr}.  After we find in the below a quantum universal $R$--matrix for
our $(A,\Delta)$, we will show that this $R$--matrix is indeed closely
related with the classical $r$--matrix given earlier (section 1, equation
\eqref{(rmatrix)}) at the Lie bialgebra level.

 Recall that the classical $r$--matrix associated with our construction
is an element in $\tilde{\Gh}\otimes\tilde{\Gh}$.  This suggests that we
better consider the Hopf $C^*$--algebra $(\tilde{A},\tilde{\Delta})$,
instead of $(A,\Delta)$.  So we need to look for our quantum $R$--matrix
in $M(\tilde{A}\otimes\tilde{A})$.  Motivated by the $R$--matrix constructed
at the QUE algebra level \cite{CGST}, \cite{ALT}, we consider $R$ as the
following (continuous) function on $\tilde{G}\times\tilde{G}$:
\begin{equation}\label{(functionR)}
R(p,q,r,s,p',q',r',s')=\bar{e}\bigl[\lambda(rs'+r's)\bigr]
\bar{e}\bigl[2\lambda(e^{-\lambda r'})p\cdot q'\bigr].
\end{equation}

 Let us try to formulate a more proper definition of $R$ as an operator on
$\tilde{\mathcal H}\otimes\tilde{\mathcal H}$.  First, let us view $R$ as
a product of two functions $\Phi$ and $\Phi'$ given by
\begin{align}
\Phi(p,q,r,s,p',q',r',s')&=\bar{e}\bigl[\lambda(rs'+r's)\bigr] \notag \\
\Phi'(p,q,r,s,p',q',r',s')&=\bar{e}\bigl[2\lambda(e^{-\lambda r'})p\cdot
q'\bigr]. \notag
\end{align}
By using partial Fourier transform purely formally and by using the 
multiplication law of $\tilde{A}$ (See Example \ref{ex} in Appendix of
section 3.), we may regard $\Phi$ and $\Phi'$ as operators on $\tilde
{\mathcal H}\otimes\tilde{\mathcal H}$:
\small 
\begin{align}
\Phi\xi(x,y,r,w,x',y',r',w')&=\xi(e^{-\lambda r'}x,e^{\lambda r'}y,r,
w-\lambda r',e^{-\lambda r}x',e^{\lambda r}y',r',w'-\lambda r) \notag \\
\Phi'\xi(x,y,r,w,x',y',r',w')&=\int\bar{e}\bigl[2\lambda(e^{-\lambda r'})
\tilde{p}\cdot\tilde{q}\bigr]e(\tilde{p}\cdot\tilde{x}+\tilde{q}\cdot
\tilde{y})\bar{e}\bigl[\eta_{\lambda}(r)\beta(\tilde{x},y)\bigr] \notag \\
&\quad\quad\quad\xi(x-\tilde{x},y,r,w,x',y'-\tilde{y},r',w')\,d\tilde{p}
d\tilde{q}d\tilde{x}d\tilde{y}.  \notag
\end{align}
\normalsize

\begin{defn}\label{R}
Define $R$ as an operator in ${\mathcal B}(\tilde{\mathcal H}\otimes\tilde
{\mathcal H})$ by $R=\Phi\Phi'$.  That is,
\small
\begin{align}
&R\xi(x,y,r,w,x',y',r',w')=\Phi\Phi'\xi(x,y,r,w,x',y',r',w')  \notag\\
&=\int\bar{e}\bigl[2\lambda(e^{-\lambda r'})
\tilde{p}\cdot\tilde{q}\bigr]e(\tilde{p}\cdot\tilde{x}+\tilde{q}\cdot
\tilde{y})\bar{e}\bigl[\eta_{\lambda}(r)\beta(\tilde{x},e^{\lambda r'}y)
\bigr]  \notag \\
&\qquad\quad\xi(e^{-\lambda r'}x-\tilde{x},e^{\lambda r'}y,r,w-\lambda r',
e^{-\lambda r}x',e^{\lambda r}y'-\tilde{y},r',w'-\lambda r)
\,d\tilde{p}d\tilde{q}d\tilde{x}d\tilde{y}.  \notag
\end{align}
\normalsize
\end{defn}

\begin{prop}\label{R1}
Let $R$ be the operator defined above.  Then $R\in M(\tilde{A}\otimes
\tilde{A})$.
\end{prop}

\begin{proof}
It is enough to show that $\Phi$ and $\Phi'$ are both left and right
multipliers.  To show this, consider an arbitrary function $F$ in the
dense subalgebra $\tilde{\mathcal A}\otimes\tilde{\mathcal A}$ of
$M(\tilde{A}\otimes\tilde{A})$, where $\tilde{\mathcal A}$ is as defined
in Example \ref{ex}.  Then by straightforward calculation, we have:
\small
\begin{align}
(\Phi F)(p,q,r,s,p',q',r',s')&=\bar{e}\bigl[\lambda(rs'+r's)\bigr]\,
F\bigl(e^{\lambda r'}p,e^{-\lambda r'}q,r,s,
e^{\lambda r}p',e^{-\lambda r}q',r',s'\bigr) \notag \\
(F\Phi)(p,q,r,s,p',q',r',s')&=\bar{e}\bigl[\lambda(rs'+r's)\bigr]\,
F(p,q,r,s,p',q',r',s'). \notag
\end{align}
\normalsize
These equations are understood to mean that $\Phi F\in{\mathcal B}
(\tilde{\mathcal H}\otimes\tilde{\mathcal H})$ is exactly the operator
realization of the function $(\Phi F)\in\tilde{\mathcal A}\otimes\tilde
{\mathcal A}$ defined by the first equation, and similarly for $F\Phi$.
From this, it is clear that $\Phi$ is a multiplier.  

 The proof that $\Phi'$ is a left multiplier follows from the following:
\small
$$
\Phi'F(p,q,r,s,p',q',r',s')=\bar{e}\bigl[2\lambda(e^{-\lambda r'})p\cdot
q'\bigr]\,F\bigl(p,q+2\lambda\eta_{\lambda}(r)e^{-\lambda r'}q',r,s,
p',q',r',s'\bigr), \notag
$$
\normalsize
which is again understood in the same way as above.  To prove that $\Phi'$
is also a right multiplier, it is more convenient to consider the Schwartz
function space in the $(p,q,r,w)$ variables having compact support both in
the $r$ and the $w$ variables, which is isomorphic (via partial Fourier
transform) to $\tilde{\mathcal A}\otimes\tilde{\mathcal A}$.  If $F$ is in
this space, we then have:
\small
\begin{align}
&F\Phi'(p,q,r,w,p',q',r',w')  \notag \\
&=\bar{e}\bigl[2\lambda(e^{-\lambda r'+w-w'})
p\cdot q')\bigr]F\bigl(p,q,r,w,p'+2\lambda(e^{-\lambda r'+w-w'})
\eta_{\lambda}(r')p,q',r',w'\bigr).  \notag
\end{align}
\normalsize
So $\Phi'$ is also a right multiplier.  Thus $R=\Phi\Phi'$ is both left
and right multiplier.
\end{proof}

\begin{prop}\label{R2}
Let $R$ be as above.  Then $R$ is an invertible element in $M(\tilde{A}
\otimes\tilde{A})$ and $R$ satisfies:
$$
R\tilde{\Delta}L_{a,b,c,d}R^{-1}=\tilde{\Delta}^{\operatorname{op}}
L_{a,b,c,d},\quad(a,b,c,d)\in\tilde{H}.
$$
Here $L_{a,b,c,d}$ denotes the ``building block'' operator defined earlier.
Thus $R$ makes $(\tilde{A},\tilde{\Delta})$ an almost cocommutative Hopf
$C^*$--algebra.
\end{prop}

\begin{proof}
It turns out that as an operator, $R$ is unitary.  And $R^*$ is:
\small
\begin{align}
&R^*\xi(x,y,r,w,x',y',r',w')  \notag\\
&=\int e\bigl[2\lambda(e^{-\lambda r'})
\tilde{p}\cdot\tilde{q}\bigr]\bar{e}(\tilde{p}\cdot\tilde{x}+\tilde{q}
\cdot\tilde{y})e\bigl[\eta_{\lambda}(r)\beta(\tilde{x},y)\bigr]  \notag \\
&\qquad\quad\xi(e^{\lambda r'}x+e^{\lambda r'}\tilde{x},e^{-\lambda r'}y,
r,w+\lambda r', e^{\lambda r}x',e^{-\lambda r}y'+e^{-\lambda r}\tilde{y},
r',w'+\lambda r), \notag
\end{align}
\normalsize
where the integration is with respect to $(\tilde{p},\tilde{q},\tilde{x},
\tilde{y})$ variables.  By using the expression for $\tilde{\Delta}
L_{a,b,c,d}$ given in Example \ref{ex}, we obtain:
\small
\begin{align}
&R\tilde{\Delta}L_{a,b,c,d}R^*\xi(x,y,r,w,x',y',r',w')  \notag\\
&=\bar{e}\bigl[\eta_{\lambda}(r)\beta(a,y-b)+\eta_{\lambda}(r')\beta
(e^{\lambda r}a, y'-e^{\lambda r}b)\bigr]\bar{e}\bigl[(r+r')c\bigr] \notag \\
&\quad\xi(e^{-d}x-e^{-d}a,e^{d}y-e^{d}b,r,w-d,e^{-d}x'-e^{\lambda r-d}a,
e^{d}y'-e^{\lambda r+d}b,r',w'-d)  \notag \\
&=\tilde{\Delta}^{\operatorname{op}}L_{a,b,c,d}
\xi(x,y,r,w,x',y',r',w').  \notag
\end{align}
\normalsize
Since the almost cocommutativity condition holds for the building blocks,
it is true for any element of $\tilde{A}$.
\end{proof}

\begin{theorem}\label{R3}
Let $R$ be defined by Definition \ref{R}.  Then $R$ satisfies the QYBE
and the quasitriangularity condition given in Definition \ref{qtrR}.
Combining this result with those of Proposition \ref{R1} and Proposition
\ref{R2}, we conclude that $R$ is a ``quasitriangular'' quantum universal
$R$--matrix for $(\tilde{A},\tilde{\Delta})$. 
\end{theorem}

\begin{proof}
The verification of the QYBE ($R_{12}R_{13}R_{23}=R_{23}R_{13}R_{12}$) is
a straightforward calculation.  We also have:
$$
\tilde{U}_{23}R_{12}\tilde{U}_{23}^*=R_{13}R_{12}\quad{\text {and}}\quad  
\tilde{U}_{12}R_{13}\tilde{U}_{12}^*=R_{13}R_{23}, 
$$
using the definition of $\tilde{U}$ given in Example \ref{ex}.  Since 
$(\operatorname{id}\otimes\tilde{\Delta})(R)=\tilde{U}_{23}R_{12}
\tilde{U}_{23}^*$ and since $(\tilde{\Delta}\otimes\operatorname{id})(R)
=\tilde{U}_{12}R_{13}\tilde{U}_{12}^*$, the quantum $R$--matrix condition
follows.  Thus we conclude that $R$ is indeed a quasitriangular quantum
universal $R$--matrix for $(\tilde{A},\tilde{\Delta})$.
\end{proof}

 Finally, let us try to relate our quantum $R$--matrix with the classical
$r$--matrix:
$$
r=\mathbf{z}\otimes\mathbf{d}+\mathbf{d}\otimes\mathbf{z}+2\sum_{i=1}^n
\mathbf{x_i}\otimes\mathbf{y_i}\in\tilde{\Gh}\otimes\tilde{\Gh}.
$$
It involves regarding $\lambda$ as a deformation parameter, and $r$ is then
a ``classical limit'' of $R$.  Since we have so far been viewing $\lambda$
as a fixed constant built into the definition of $G$ and its Poisson bracket,
we have to approach a little differently.  It actually corresponds to the
deformation process of the dual Hopf $C^*$--algebra $\hat{A}$.

 One serious problem is that as we try to let $\lambda$ vary, the algebra
$\tilde{A}$ (or $\tilde{A}\otimes\tilde{A}$) also changes.  Because of this,
we will only work on its dense function space $\tilde{\mathcal A}$ (or
$\tilde{\mathcal A}\otimes\tilde{\mathcal A}$), ignoring its algebra structure.
Again, as in the proof of Proposition \ref{R1}, it is more convenient to
regard $\tilde{\mathcal A}$ as the functions in the $(p,q,r,w)$ variables:
That is, the Schwartz function space having compact supportin the $r$ and
the $w$ variables.  The $L^1$--completion of $\tilde{\mathcal A}\otimes
\tilde{\mathcal A}$ is isomorphic to $L^1(\tilde{H}\otimes\tilde{H})$,
independent of the value of $\lambda$.

 Recall that we could realize $R$ as a continuous function on $\tilde{G}
\times\tilde{G}$ by equation \eqref{(functionR)}.  To emphasize its dependence
on $\lambda$, let us denote it from now on by $R_{\lambda}$.  Consider the
operator $\Psi_{\lambda}$ on $\tilde{\mathcal A}\otimes\tilde{\mathcal A}$
(for the time being, $\tilde{\mathcal A}$ is viewed as an algebra) defined by
$$
\Psi_{\lambda}(F)=R_{\lambda}FR_{\lambda}^*,
\qquad F\in\tilde{\mathcal A}\otimes\tilde{\mathcal A}.
$$
Then we have:
\small
\begin{align}
&\Psi_{\lambda}(F)(p,q,r,w,p',q',r',w')=\bar{e}\bigl[2\lambda
(e^{-\lambda r})p\cdot q'\bigr]e\bigl[2\lambda(e^{w-w'-\lambda r})p\cdot q'
\bigr]  \notag \\
&\quad F\bigl(e^{\lambda r'}p,e^{-\lambda r'}q+2\lambda e^{-\lambda r
-\lambda r'}\eta_{\lambda}(r)q',r,w,e^{\lambda r}p'-2\lambda e^{w-w'}
\eta_{\lambda}(r')p,e^{-\lambda r}q',r',w'\bigr).   \label{(PsiF)}
\end{align}
\normalsize
By $L^1$--extension, we will define $\Psi_{\lambda}$ as an operator on the
Banach space $L^1(\tilde{H}\otimes\tilde{H})$, ignoring any algebra structure,
via the equation \eqref{(PsiF)}.  This would be our operator realization of
$R_{\lambda}$.

 Let us now consider the classical $r$--matrix.  First, by means of the dual
pairing between $\tilde{\Gh}^*$ and $\tilde{\Gh}$, we may regard $r\in\tilde
{\Gh}\otimes\tilde{\Gh}$ as a linear function on $\tilde{\Gh}^*\otimes\tilde
{\Gh}^*$.  Let us denote it by $\psi$:
$$
\psi(p,q,r,s,p',q',r',s')=rs'+r's+2p\cdot q'.
$$
Next, we have to find a way to make $\psi$ to determine an operator on $L^1(
\tilde{H}\otimes\tilde{H})$.  Since it should correspond to $\lambda=0$ case,
we will construct an (unbounded) operator such that it looks like an (unbounded)
``derivation'' with respect to the multiplication (for $\lambda=0$) on $\tilde
{\mathcal A}\otimes\tilde{\mathcal A}$.  That is, we consider the densely
defined operator:
$$
F\mapsto [\psi,F]=\psi\times_{\lambda=0}F-F\times_{\lambda=0}\psi,
\qquad F\in\tilde{\mathcal A}\otimes\tilde{\mathcal A}.
$$
But $\times_{\lambda=0}$ is essentially the ordinary convolution on $S(\tilde
{H})$ (or $S(\tilde{H}\times\tilde{H})$).  So by straightforward calculation,
again formally using partial Fourier transform, we obtain:
\small
\begin{align}
&[\psi,F](p,q,r,w,p',q',r',w') \notag \\
&=\int\bigl[(r\tilde{\tilde{s}}+r'\tilde{s})+2(p\cdot q'
+rq'\cdot\tilde{y}-e^{w-w'}p\cdot q'-r'e^{w-w'}p\cdot\tilde{x})
\bigr]  \notag\\
&\quad\quad e[\tilde{s}\tilde{w}+\tilde{\tilde{s}}\tilde{\tilde{w}}]
e[\tilde{p}\cdot\tilde{x}+\tilde{q}\cdot\tilde{y}]F(e^{\tilde{w}}p,
e^{-\tilde{w}}q+\tilde{q},r,w,e^{\tilde{\tilde{w}}}p'+\tilde{p},
e^{-\tilde{\tilde{w}}}q',r',w'),  \label{[psiF]}
\end{align}
\normalsize
where the integration is taken with respect to all the tilde ($\tilde{\ }$)
and double tilde ($\tilde{\tilde{\ }}$) variables.  From now on, we will just
use \eqref{[psiF]} as our defining equation for $[\psi,\cdot]$, an unbounded
operator on the Banach space $L^1(\tilde{H}\otimes\tilde{H})$.  This would be
our operator realization of $\psi$.

 Then by comparing the formulas \eqref{(PsiF)} and \eqref{[psiF]}, we obtain
the following result.  Although we showed directly in Theorem \ref{R3} that
our $R$ satisfies the QYBE, this proposition indicates that this property is
actually suggested by the CYBE satisfied by the associated classical $r$-matrix.
 
\begin{prop}
Let the the notation be as above.  Then: 
$$
\lim_{\lambda\to0}\left\|\frac{\Psi_{\lambda}(F)-F}{\lambda}-(-2\pi i)
[\psi,F]\right\|_{L^1}=0,
$$
for $F\in\tilde{\mathcal A}\otimes\tilde{\mathcal A}$.  Thus at least
in the sense of the operators on the Banach space $L^1(\tilde{H}\times
\tilde{H})$, we may say that the ``classical limit'' as $\lambda\to0$
of our quantum $R$--matrix $R_{\lambda}$ is $(-2\pi i)\psi$, the operator
realization of the classical $r$--matrix.
\end{prop}

\begin{proof}
From equation \eqref{(PsiF)}, we may express $\Psi_{\lambda}(F)$ as
follows, taking advantage of the Fourier inversion theorem:
\small
\begin{align}
&\Psi_{\lambda}(F)(p,q,r,w,p',q',r',w') \notag \\
&=\int\bar{e}\bigl[\lambda(r\tilde{\tilde{s}}+r'\tilde{s})\bigr]
\bar{e}\bigl[2\lambda(e^{-\lambda r})p\cdot q'\bigr]\bar{e}\bigl[
2\lambda(e^{-\lambda r-\lambda r'})\eta_{\lambda}(r)q'\cdot\tilde{y}
\bigr]  \notag\\
&\quad\quad e\bigl[2\lambda(e^{w-w'-\lambda r})p\cdot q'\bigr]
e\bigl[2\lambda(e^{w-w'})\eta_{\lambda}(r')p\cdot\tilde{x}\bigr]
e[\tilde{s}\tilde{w}+\tilde{\tilde{s}}\tilde{\tilde{w}}]
e[\tilde{p}\cdot\tilde{x}+\tilde{q}\cdot\tilde{y}]  \notag\\
&\quad\quad F(e^{\tilde{w}}p,e^{-\tilde{w}}q+\tilde{q},r,w,
e^{\tilde{\tilde{w}}}p'+\tilde{p},e^{-\tilde{\tilde{w}}}q',r',w'). \notag
\end{align}
\normalsize
The integration is with respect to all the tilde and double tilde variables.
Comparing this expression with equation \eqref{[psiF]} for $[\psi,F]$, we can
see easily the pointwise convergence.  The $L^1$ convergence is proved using
the Lebesgue's dominated convergence theorem.
\end{proof}

 The quantum universal $R$--matrix is useful in the study of representation
theory of our Hopf $C^*$--algebras $(\tilde{A},\tilde{\Delta})$ and
$(A,\Delta)$.  We will study representation theory of our quantum groups
elsewhere (See \cite{BJKpz}.).  It turns out that the representation theory
satisfies interesting {\em quasitriangularity\/} property, which is not
present in the earlier examples of quantum groups corresponding to linear
Poisson brackets.


\bibliography{ref}

\bibliographystyle{amsplain}

\end{document}